\documentclass[11pt,a4paper,reqno]{amsart}
\usepackage{a4wide}
\usepackage{amsaddr}
\usepackage{amsmath}
\usepackage{amsfonts, amssymb, amsthm}
\usepackage{paralist}
\usepackage{enumitem}
\usepackage[colorlinks=true]{hyperref}
\hypersetup{urlcolor=blue, citecolor=red}

\makeatletter
\renewcommand{\email}[2][]{%
	\ifx\emails\@empty\relax\else{\g@addto@macro\emails{,\space}}\fi%
	\@ifnotempty{#1}{\g@addto@macro\emails{\textrm{(#1)}\space}}%
	\g@addto@macro\emails{#2}%
}
\makeatother

\theoremstyle{definition}

\numberwithin{equation}{section}

\usepackage[scrtime]{prelim2e}

\usepackage[normalem]{ulem}
\usepackage{color}
\begin{document}
\title[The Choquard equation with mass supercritical nonlinearity]
{Normalized solutions for the Choquard equation with mass supercritical nonlinearity}
\author{Na Xu}\email{weila319@163.com}
\address{School of Science,
Tianjin University of Technology and Education\\
    Tianjin 300222, China}
\author{Shiwang Ma$^*$}\email{shiwangm@nankai.edu.cn}
\address{School of Mathematical Sciences and LPMC, Nankai University\\
	Tianjin 300071, China}


\thanks{*Corresponding author: Shiwang Ma\\
Email Addresses: shiwangm@nankai.edu.cn (SM)\\
}



\date{}

\begin{abstract}
 We consider the  nonlinear Choquard equation
$$\begin{cases}
& - \Delta  u   = (I_\alpha \ast F(u))F'(u) -\mu u \  \text{in}\  \mathbb{R}^N,\\
&  u   \in \ H^1(\mathbb{R}^N), \  \int_{\mathbb{R}^N} |u|^2 dx=m,
\end{cases}  $$
where  $\alpha\in(0,N)$, $m>0$ is prescribed, $\mu \in \mathbb{R}$ is a Lagarange multiplier,
and $I_\alpha$ is the Riesz potential.
 Under some mild mass supercritical and Sobolev subcritical  conditions,    we prove the existence and multiplicity of normalized solutions.
 \end{abstract}

\maketitle


\section{Introduction and main results}

In this paper, we consider the  Choquard equation
$$\begin{cases}
& - \Delta  u   = (I_\alpha \ast F(u))f(u) -\mu u \  \quad  \text{in}\  \mathbb{R}^N,\\
&   u   \in \ H^1(\mathbb{R}^N), \   \int_{\mathbb{R}^N} |u|^2 dx=m.\\
\end{cases} \eqno(1.1)$$
Here $N\geq 3,$  $\alpha\in(0,N)$, $f\in C(\mathbb{R},\mathbb{R}),$ $F(t)=\int_0^t f(s)ds$,  $m>0$ is a given constant,
  $\mu \in \mathbb{R}$ will arise as a Lagarange multiplier,
  and $I_\alpha:\mathbb{R}^N\rightarrow \mathbb{R} $ is the Riesz potential defined for every $x\in \mathbb{R}^N \setminus \{0\}$ by
$$ I_\alpha(x)=\frac{\Gamma(\frac{N-\alpha}{2})}{ \Gamma(\frac{\alpha}{2}) \pi^{N/2}2^\alpha | x|^{N-\alpha}},
 $$
 where  $\Gamma$ is the Gamma function  see \cite{Riesz}.

The Choquard equation
$$
 - \Delta  u +u  = (I_2 \ast |u|^2)u  \  \text{in}\  \mathbb{R}^3, \    \eqno(1.2)$$
has been introduced  by S. I. Pekar in 1954 \cite{Pekar} as  a physical
model describing  the quantum mechanics of a polaron at rest.
In \cite{Lieb}, P. Choquard applied it as an approximation to Hartree-Fock theory of one component plasma.
Moreover, \cite{Gross} and \cite{Penrose} used (1.2) in  multiple particle systems  and quantum mechanics.
 Indeed, if $u$ is a solution of (1.2), then the wave
function $\psi(x,t) = e^{-it}u(x)$ is a solitary one of the focusing time-dependent Hartree equation
$$
i\frac{\partial}{\partial t}\psi +\Delta \psi +(I_2\ast |\psi|^2)\psi =0, \quad (x,t)\in \mathbb R^3\times \mathbb R.
$$

 E. H. Lieb, P. L. Lions and G. Menzala \cite{Lieb,Lions,Menzala} studied  problem (1.2)   by establishing  variational framework. Since then, there have been many papers considering problem (1.2) or a similar problem with general pure power nonlinearity
$$
-\Delta u + u = (I_\alpha \ast |u|^p)|u|^{p-2}u, \quad x\in \mathbb R^N,
$$
where $N \ge  3$ and $\alpha\in (0,N)$ by variational methods. Until now the existence results have been mostly available when the nonlinearity $F(u)$ is homogeneous.
 In recent decades, the problem of finding normalized solutions of nonlinear Choquard type  equations
  has received considerable attention,  such as \cite{Genev,Feng,Cazenave2,Ye,Li,Li2,Liu,Yao,Li3,Yao2,He}  and the references therein.

Define the energy functional $I: H^1(\mathbb{R}^N) \rightarrow \mathbb{R}$ by
$$
I(u)=\frac{1}{2}\int_{\mathbb{R}^N} | \nabla u|^2dx-\frac{1}{2}\int_{\mathbb{R}^N} (I_\alpha  \ast F(u))F(u)dx.
\eqno(1.3)
$$
Then for $u \in H^1(\mathbb{R}^N),$ we have
$$ I'(u) \varphi
= \int_{\mathbb{R}^N}  \nabla u  \nabla \varphi dx- \int_{\mathbb{R}^N} (I_\alpha  \ast F(u))f(u) \varphi dx ,
\ \text{for all } \  \varphi \in C_0^\infty (\mathbb{R}^N)$$
and a critical point of $I$ constrained
to
$$ S_m= \left\{u\in H^1(\mathbb{R}^N):\int_ {\mathbb{R}^N} |u|^2dx=m\right\}$$
gives rise to a solution to (1.1).
The  variational methods  are heavily dependent on the behavior of the nonlinearity when seeking for normalized solutions to the Choquard equation.
In the mass subcritical case, the constrained functional $I|_{S_m}$ is bounded from below and coercive.
 Under general mass subcritical  conditions of Berestycki-Lions type, Cingolani and Tanaka \cite{Cingolani2} obtained the existence of ground states of (1.1) and
 \cite{Cingolani1}  studied the  existence of infinitely many   normalized solution solutions for (1.1) with  $L^2$  subcritical growth  at $\infty$ by using minimax methods.

 However, in  the mass supercritical case,  $I$ is unbounded
from below on $S_m$ for any $m > 0.$
Therefore, more difficulties lie ahead and few results are obtained  in the mass supercritical case.
  Li and Ye \cite{Li5} considered the problem (1.1)  under the following conditions:\\
 $(f_0)$ $f \in C^1(\mathbb{R},\mathbb{R} )$ and  $f(s)\equiv0$ for $s\leq 0$;  \\
 $(f_1)$ there exists $r\in(\frac{N+\alpha +2}{N},\frac{N+\alpha  }{N-2})$ such that
$$ \lim_{|s| \rightarrow 0} \frac{f(s)}{|s|^{r-2}s}=0 \quad \text{and} \ \quad \lim_{|s| \rightarrow +\infty} \frac{F(s)}{|s|^{r } }=+\infty;$$
$(f_2)$ $\lim_{|s| \rightarrow +\infty} \frac{F(s)}{|s|^{\frac{N+\alpha  }{N-2}} }=0$;\\
$(f_3)$ there exists $\theta_1\geq 1$ such that $\theta_1 \hat{F}(s) \geq \hat {F}(ts) $ for $s\in \mathbb{R}$ and $t\in[0,1],$
where $$ \hat{F}(s) =f(s)s-\frac{N+\alpha +2}{N}F(s);
$$
$(f_4)$ $f(s)s< \frac{N+\alpha  }{N-2}F(s)$ for all $s >0;$\\
$(f_5)$ $\tilde {F}'(s) $ exists and
$\tilde {F}'(s)  s >  \frac{N+\alpha +2}{N}\tilde {F}(s)$, where
 $$\tilde {F}(s):= f(s)s-\frac{N+\alpha  }{N}F(s).\eqno(1.4)
 $$
$(f_6)$ there exists $0<\theta_2<1$ and $t_0>0 $ such that for all  $s\in \mathbb{R}$ and $|t| \leq t_0,$
$$F(ts)\leq \theta_2 |t|^{\frac{N+\alpha +2 }{N}}F(s).$$
By using the methods in \cite{Jeanjean2} and  the concentration
compactness due to P. L. Lions,  Li and Ye   obtained  the existence of positive normalized solutions for (1.1).
Subsequently,  \cite{Yuan} improved their results and  obtained the existence of positive normalized solutions without the condition $f\in C^1(\mathbb{R},\mathbb{R})$ by using a minimax procedure established by Jeanjean \cite{Jeanjean2} and Chen and Tang \cite{Chen2}.

By using the mountain pass theorem,   Bartsch et al. \cite{Bartsch2} obtained the existence  of  ground sate normalized solutions to (1.1) under the following conditions:\\
$(g_0)$ $f\in C^0(\mathbb{R},\mathbb R);$ \\
$(g_1)$  there exists $r,p \in \mathbb{R}$ verifying
$\frac{N+\alpha+2 }{N}<r \leq p<\frac{N+\alpha}{N-2}$
such  that
$$0<rF(s)\leq f(s)s \leq pF(s) \quad \text{for} \  s\neq0;$$
$(g_2)$
$ \tilde {F}(s)/|s|^{ 1+\frac{\alpha+2}{N}}$ is   nonincreasing on $(-\infty,0)$
and nondecreasing  on $(0,\infty)$, where
$\tilde {F}(t)$ is given in (1.4).\\
If, in addition, $f$ is odd, the authors in \cite{Bartsch2} also obtained an unbounded sequence of pairs of radial normalized solutions by using the linking theorem.

In  \cite{Jeanjean4}, Jeanjean and Lu concerned  with the following the nonlinear scalar field
equation with  $L^2$ constraint
$$\begin{cases}
& - \Delta  u   = g(u) -\mu u \  \quad  \text{in}\  \mathbb{R}^N,\\
&   u   \in \ H^1(\mathbb{R}^N), \   \int_{\mathbb{R}^N} |u|^2 dx=m,\\
\end{cases}
\eqno(1.5)$$
 Assuming only that the nonlinearity $g$ is
continuous and satisfies weak mass supercritical conditions, Jeanjean and Lu obtained the existence of ground
states to (1.5) and reveal the basic behavior of the ground state energy $E_m$ as $m >0$ varies.
Moreover, they also obtain infinitely
many radial solutions for any $N\ge 2$ and established the existence and multiplicity of nonradial
sign-changing solutions for $N\ge 4$.

In the present paper,  inspired by \cite{Jeanjean4}, we study the existence and multiplicity of normalized solutions of (1.1) under weaker $L^2$-supercritical conditions.
To this end,   we make  the following assumptions.\\
${\bf(H0)}$ ($continuity  \ condition$)  $f :\mathbb{R} \rightarrow \mathbb{R} $  is continuous.\\
${\bf (H1)}$ ($L^2$  $supercritical \ condition$)  $   \lim _{t\rightarrow 0}f(t)/|t|^{ \frac{\alpha+2}{N}}=0$ and  $   \lim _{t\rightarrow \infty}F(t)/|t|^{ 1+\frac{\alpha+2}{N}}=+\infty. $\\
${\bf (H2)}$  ($Sobolev \ subcritical \ condition \ at \ \infty$) $   \lim_{|t| \rightarrow +\infty} F(t)/|t|^{\frac{N+\alpha  }{N-2}}=0 .$\\
${\bf (H3)}$  ($global \ condition$)  $ f( t)t< \frac{N+\alpha }{N-2}F(t)  $ for all $t\in  \mathbb{R} \backslash \{0\}$ and the map $t\mapsto \tilde{F}(t)/|t|^{ 1+\frac{\alpha+2}{N}} $ is   strictly decreasing on $(-\infty,0)$
and strictly increasing on $(0,\infty),$
where $\tilde{F}(t)$ is given in (1.4). \\
${\bf (H4)}$ ($Sobolev \ subcritical \ condition \ at \ 0$)  $\lim _{t\rightarrow 0}F(t)/|t|^{ \frac{N+\alpha }{N-2}}=+\infty. $

The assumptions $(H0)-(H4)$ are similar to that used in \cite{Jeanjean4}. Plainly,  the $L^2$ supercritical condition $(H1)$ is weaker than the condition $(f_1)$ in \cite{Li5}.
It is easy to see that the monotonicity condition in $(H3)$ is weaker than $(f_5)$ and hence the global condition $(H3)$ sharply improves the global conditions $(f_3)-(f_6)$ in \cite{Li5}.
The $L^2$ supercritical condition $(H1)$ and Sobolev subcritical condition  $(H2)$ are also weaker than the corresponding one in \cite{Bartsch2}.
 Moreover, since the absence of $(f_3)$ in our assumptions,  some  new
difficulties arise and the technique used  in \cite{Li5,Yuan} cannot be used any more.
In  this paper,  we  mainly   adopt  the minimax theorem which  is different   from  that  used in \cite{Bartsch2,Li5,Yuan} and extend some results in \cite{Jeanjean4} to (1.1).  Due to the nonlocal feature of our problem (1.1), some careful analysis is needed.

As given in \cite{Jeanjean4}, an explicit example can be constructed as follows.
Setting
$ \alpha _N:=\frac{4+\alpha }{N(N-2)}, $
and define the odd continuous function
$$ f(t):=\left[(1+\frac{\alpha+2}{N})\ln(1+ |t|^{\alpha_N})+ \frac{\alpha_N |t|^{\alpha_N}}{1+ |t|^{\alpha_N}}\right]|t|^{\frac{(\alpha+2-N)}{N}}t,$$
with the primitive function $F(t):= |t|^{1+\frac{ \alpha+2 }{N}}\ln(1+ |t|^{\alpha_N}).$
We can see  that the function $f,F$ satisfy the conditions $(H0)-(H4) $  but  not  the  ones  in  \cite{Bartsch2}.
Moreover, noting that $ \lim_{|s|\rightarrow \infty}\frac{F(s)}{|s|^r}=0,$ the function $F$ does not satisfy the conditions $( f_1).$

For any given $m>0,$  define
$$
E_m=\inf_{u \in \mathcal{P}_m} I(u),
\eqno(1.6)
$$
where $ \mathcal{P}_m$ is the Poho\v zaev manifold defined by $\mathcal{P}_m :=\left\{  u \in S_m: \ P(u)=0 \right\}$ and
$$
P(u):=
 \int_{\mathbb{R}^N} | \nabla u|^2dx -\frac{N}{2}\int_{\mathbb{R}^N} (I_\alpha  \ast F(u))\left(f(u)u-\frac{N+\alpha}{N}F(u)\right)dx.
\eqno(1.7)
$$

Our main results are as follows.

{\bf Theorem 1.1.} {\it Let  $f$ satisfy $(H0)-(H3).$
 Then for any $m>0$, $(1.1) $ admits a ground state and the associated Lagrange multiplier $\mu$ is positive.
If  $f$  is odd, then the ground state is positive.
Furthermore, in addition to $(H0)-(H3)$,  if $(H4)$ also holds, then the least energy $E_m$ is positive,
continuous, strictly decreasing with respect to $m$, and satisfies
$$\lim_{m \rightarrow 0^+} E_m = +\infty  \ \text{and}  \  \lim_{m \rightarrow +\infty} E_m =0.$$ }

{\bf Theorem 1.2.} {\it
Assume $f$ is odd and $(H0)-(H3)$ hold.  Then $(1.1)$ has infinitely
many radial solutions $\{u_k \}_{k=1}^\infty$ for any $m > 0.$
In particular,
$$I (u_{k+1}) \geq I (u_k) > 0   \  \   \text{for each } k\in \mathbb{N}^+, $$
  and $I (u_k )\rightarrow +\infty$   as $k \rightarrow \infty.$}

 \smallskip

\textbf{Organization of the paper}.  In Section 2, we give some preliminaries.
 In Section 3, we consider the properties of $E_m.$
 Section 4  devoted to the proofs of Theorems 1.1  and 1.2.
In Section 5, we  study the multiplicity   result of radial solution.

\smallskip

\textbf{Basic notations}. Throughout this paper, we assume $N\geq
3$. \begin{itemize}
   \item For $p\in [1,+\infty),$ let $L^p(\mathbb{R}^N)$ be the usual Lebesgue space endowed with the norm
 $\|u\|_p=\left(\int_{\mathbb{R}^N}|u|^pdx\right)^{\frac{1}{p}}.$
   \item $ H^1(\mathbb{R}^N)$  is the usual Sobolev space endowed with the norm
$  \|u\|_{H^1}^2= \|\nabla u\|_2^2+\|  u\|_2^2.$
$D^{1,2}(\mathbb{R}^N)=\left\{ u\in L^{2^*}(\mathbb{R}^N)| \   \int_{\mathbb{R}^N}| \nabla u|^2dx<+\infty\right\}.$
And $ H_r^1(\mathbb{R}^N)$ denotes the subspace of functions in $ H^1(\mathbb{R}^N)$
 which are radially symmetric with respect to zero.
   \item For any $m>0,$ let
$$ S_m:=\left\{u\in H^1(\mathbb{R}^N) | \   \|u\|_{2}^2=m\right\},
\ \text{and}\
  B_m:=\left\{ u\in H^1(\mathbb{R}^N) |  \  \|u\|_{2}^2 \leq m\right\}.$$
    \item $C,$ $ C_1,$  $ C_2,$ $\cdots$ denote positive constants, whose values can change from line to line.
 \end{itemize}

\section{Preliminary results}

We first recall the Hardy-Littlewood-Sobolev inequality \cite{Lieb2}.

\quad {\bf Lemma 2.1.} {\it
If  $1<r,t< \infty ,$
 and $\alpha\in(0, N)$  with $\frac{1}{r}+\frac{1}{t}=1 +\frac{\alpha}{N}. $
For $ \varphi \in L^r(\mathbb{R}^N)$ and $\psi \in L^ t(\mathbb{R}^N), $
then there exists a sharp constant $C(N, \alpha, r, t) >0$ such that
$$\left | \int_{\mathbb{R}^N}\int_{\mathbb{R}^N} \frac{\varphi(x)  \psi(y) }{ |x-y |^{N-\alpha }}dx dy\right|
\leq  C(N, \alpha, r, t)  \|\varphi\|_{r} \|\psi  \|_{ t}.
\eqno(2.1)
$$
}
Next, we present some technique lemmas which are needed in the proofs of our main results.

\quad{\bf Lemma 2.2.} {\it  Let $(H0)-(H3)$  be satisfied. Then the following statements
hold.\\
(i ) For any $m > 0,$ there exists $\delta =\delta(N,m) > 0 $ small enough such that
$$  \frac{1}{4}\int_{\mathbb{R}^N} | \nabla u|^2dx
 \leq I(u)\leq  \frac{1}{2}\int_{\mathbb{R}^N} | \nabla u|^2dx$$
 for all $u\in B_m$ satisfying $ \| \nabla u \|_2\leq \delta.$\\
(ii ) Let $ \{u_n\}$ be a bounded sequence in $H^1(\mathbb{R}^N)  .$
 If $u_n \rightarrow 0$ in $L^{\frac{2(N+\alpha+2)}{N+\alpha}}(\mathbb{R}^N) ,$
   then
 $$ \lim_{n \rightarrow \infty}    \int_{\mathbb{R}^N} (I_\alpha \ast F(u_n))F(u_n)dx= \lim_{n \rightarrow \infty}    \int_{\mathbb{R}^N} (I_\alpha \ast F(u_n))\widetilde{F}(u_n)dx=0.$$
(iii ) Let $\{u_n\},$ $ \{v_n\}$  be bounded sequences in $H^1(\mathbb{R}^N) .$
If $v_n \rightarrow 0$ in $L^{\frac{2(N+\alpha+2)}{N+\alpha}}(\mathbb{R}^N) ,$
 then
$$   \lim_{n \rightarrow \infty}   \int_{\mathbb{R}^N} (I_\alpha \ast F(u_n))f(u_n)v_n dx=0.  $$  }
\qquad{\bf Proof.}
(i ) We show that there exists $\delta = \delta(N,m) > 0$  small enough such that
$$   \int_{\mathbb{R}^N }(I_\alpha \ast F(u ))F(u )dx \leq \frac{1}{4} \int_{\mathbb{R}^N} | \nabla u|^2dx
 \  \text{for any}  \   u \in  B_m \   \text{ with}   \   \| \nabla  u\|_{2}  \leq \delta   .
 \eqno(2.2)$$
 By $(H0)-(H2)$,  for every $ \varepsilon > 0,$ there exists
$C_\varepsilon  > 0 $ such that $|F(t)| \leq \varepsilon |t|^{1+ \frac{\alpha+2}{N}}
+ C_\varepsilon  |t|^{\frac{N+\alpha}{N-2}}$   for all $t \in \mathbb{R}. $
 In view of   Gagliardo-Nirenberg inequality, for any $ u \in  B_m ,$ one then has
\begin{align*}
 \int_{\mathbb{R}^N }| (I_\alpha \ast F(u ))F(u ) | dx
 \leq & C_1  \left(\int_{\mathbb{R}^N } | F(u )  |^{\frac{2N}{N+\alpha }}  dx\right)^{\frac{N+\alpha}{N }}\\
 \leq & C_2 \varepsilon^{2} \left(\int_{\mathbb{R}^N } | u |^{\frac{2(N+\alpha+2) }{N+\alpha}}dx\right)^{\frac{N+\alpha}{N }}
 + C_2 C_\varepsilon^{2} \left(  \int_{\mathbb{R}^N } | u |^{\frac{2 N  }{N-2 }}   dx\right)^{\frac{N+\alpha}{N }} \\
 \leq & C_3 C_{N} \varepsilon^{2}\| \nabla u \|_{2}^{2 }
 +  C_3 C_\varepsilon^{2}  C_N' \| \nabla u \|_{ 2}^{{\frac{2(N+\alpha)}{N-2}}}\\
 \leq & \left(C_4  C_{N}\varepsilon^{2}  +  C_4 C_\varepsilon^{2}  C_N' \| \nabla u \|_{2}^{\frac{2( \alpha +2)}{N-2 }} dx \right)    \| \nabla u \|_{ 2}^{2}.
\end{align*}
Choose
$$ \varepsilon:= \sqrt{\frac{1}{8C_4C_{N}}} ,  \   \delta =\left(\frac{1}{8 C_4 C_\varepsilon^{2}  C_N'}\right)^{\frac{N-2}{2(\alpha+2)}}.$$
(2.2) holds  and then Item (i) follows.

(ii)
The assumptions  $(H0)-(H2)$ imply that  for any $ \varepsilon > 0,$ there exists $C_\varepsilon>0$ such that
 $  |F(t)| \leq C_\varepsilon |t|^{1+ \frac{\alpha+2}{N}}+\varepsilon |t|^{{\frac{N+\alpha}{N-2}}}$
  for all $t \in \mathbb{R}.$
  Thus, by Lemma 2.1, we have
   \begin{align*}
  \left|  \int_{\mathbb{R}^N} (I_\alpha \ast F (u_n))\widetilde{F}(u_n)dx\right|
   \leq & C_1\left( \int_{\mathbb{R}^N} | F (u_n)|^{\frac{2N}{N+\alpha} }dx\right)^{ \frac{N+\alpha}{2N}}
   \left( \int_{\mathbb{R}^N} |\widetilde{F}(u_n)|^{\frac{2N}{N+\alpha} }dx\right)^{ \frac{N+\alpha}{2N}}  \\
    \leq & C_2 C_\varepsilon^2 \|  u_n\|_{  \frac{2(N+\alpha+2)}{N+\alpha}}^ \frac{2(N+\alpha+2)}{N}+C_2 \varepsilon^2  \|  u_n \|_{\frac{2 N }{N-2}}^{ \frac{2(N+\alpha )}{N-2}}\\
    \rightarrow &0.
\end{align*}
The proof of  $ \lim_{n \rightarrow \infty}    \int_{\mathbb{R}^N} (I_\alpha \ast F(u_n))F(u_n)dx= 0 $ is similar.

(iii)
For any $\varepsilon>0,$ by  $(H0)-(H2),$   there exists $C_\varepsilon>0$ such that
$|f(t)| \leq C_\varepsilon |t|^{{(\alpha+2)/N}}+ \varepsilon  |t|^{{(\alpha+2)/(N-2)}}$ for all $t\in \mathbb{R}.$
Hence,  we have
\begin{align*}
 \left| \int_{\mathbb{R}^N} (I_\alpha \ast  F (u_n)) f (u_n)v_ndx\right|
    \leq &   C_1\left( \int_{\mathbb{R}^N} | F(u_n) |^{\frac{2N}{N+\alpha} }dx\right)^{  \frac{N+\alpha}{2N} }\left( \int_{\mathbb{R}^N} | f(u_n)v_n|^{\frac{2N}{N+\alpha} }dx\right)^{  \frac{N+\alpha}{2N} }\\
    \leq & [ C_2  C_\varepsilon \|  u_n\|_{  \frac{2(N+\alpha+2)}{N+\alpha}}^{\frac{N+\alpha+2}{N }}
    +C_2 \varepsilon  \|  u_n \|_{\frac{2 N  }{  N-2}} ^ \frac{N+\alpha}{N-2} ]\\
   &  \cdot[ C_2  C_\varepsilon \|  u_n\|_{  \frac{2(N+\alpha+2)}{N+\alpha}}^{\frac{\alpha+2}{N}}  \|  v_n\|_{  \frac{2(N+\alpha+2)}{N+\alpha}}
    +C_2 \varepsilon  \|  u_n \|_{\frac{2 N  }{  N-2}} ^ {\frac{\alpha+2}{N-2}}
   \|  v_n \|_{\frac{2 N}{N-2}} ]\\
    \rightarrow &0,
\end{align*}
which yields the conclusion. \hfill$\Box$

For every $u \neq 0$ and $s\in \mathbb{R},$ we  introduce the map $(s \star u)(x)=e^{Ns/2}u(e^s x).$

{\bf Lemma 2.3.} {\it  Let $(H0)-(H3)$  be satisfied. Then  for any $u\in H^1(\mathbb{R}^N ) \backslash \{0\},$ one
has\\
(i ) $I (s\star u) \rightarrow 0^+ $ as $s\rightarrow -\infty,$\\
(ii ) $I (s\star u)\rightarrow -\infty$  as $s\rightarrow +\infty.$}

{\bf Proof.}
Let $m:= \|u\|_2^2>0.$
Note that $s \star u \in S_m \subset B_m$ and
$$   \|\nabla (s \star u)  \|_2=e^s \|  \nabla u\|_2 .$$
Applying Lemma 2.2(i), we obtain
$$  \frac{1} {4}e^{2s}\|  \nabla u\|_2^2 \leq I(s \star u)\leq  \frac{1} {2} e^{2s}\|  \nabla u\|_2^2 \ \text{when} \ s\rightarrow -\infty.   $$
Thus, $ \lim_{s \rightarrow -\infty} I(s \star u)=0^+.$

(ii) For every $ \lambda  \geq 0,$
   define a function $h_\lambda  : \mathbb{R}  \rightarrow \mathbb{R} $  by
$$ h_\lambda(t):=
 \begin{cases} \frac{F(t)}{|t|^{1+\frac{\alpha+2}{N}}}+\lambda , \  & \text{for} \ t \neq 0,\\
  \lambda,\qquad  & \text{for} \ t=0.\end{cases} \eqno(2.3)$$
It is easy to check that $F(t) =h_\lambda (t)|t|^{1+\frac{\alpha+2}{N}} - \lambda|t|^{1+\frac{\alpha+2}{N}}$
for all $t \in \mathbb{R}  .$
 Also,
  $(H0)-(H1)$  imply that  $h_\lambda $ is continuous and
$$    h_\lambda(t) \rightarrow  +\infty \  \text{as} \  t \rightarrow \infty. $$
Take $ \lambda > 0$ large enough satisfying  $h_\lambda(t)\geq0$ for any  $t \in \mathbb{R}  .$
Applying Fatou's lemma, we have
$$ \lim_{s \rightarrow +\infty} \int_{ \mathbb{R}^N } (I_\alpha \ast [h_\lambda(e^{Ns/2 }u)|u|^{1+\frac{\alpha+2}{N}}] )(\frac{h_\lambda(e^{Ns/2 }u)}{2}-\lambda)|u|^{1+\frac{\alpha+2}{N}} =+\infty, \eqno(2.4)$$
and
$$ \lim_{s \rightarrow +\infty} \int_{ \mathbb{R}^N } (I_\alpha \ast [(\frac{h_\lambda(e^{Ns/2 }u)}{2}- \lambda)|u|^{1+\frac{\alpha+2}{N}}] )( h_\lambda(e^{Ns/2 }u)  )|u|^{1+\frac{\alpha+2}{N}} =+\infty.\eqno(2.5)$$
Thus,  (2.4) and (2.5)  give that
\begin{align*}
  I(s \star u)
  = & \frac{1}{2}e^{2s}\left[\int_{\mathbb{R}^N} | \nabla  u |^2dx - \lambda^2   \int_{\mathbb{R}^N} (I_\alpha\ast |   u |^{1+\frac{\alpha+2}{N}})  |   u |^{1+\frac{\alpha+2}{N}} dx
\right. \\
   &\left.-\int_{\mathbb{R}^N} \int_{\mathbb{R}^N}  \frac{h_\lambda(e^{Ns/2}u(x))|u(x)|^{1+\frac{\alpha+2}{N}} (\frac{h_\lambda(e^{Ns/2} u(y))}{2}- \lambda)| u(y)|^{1+\frac{\alpha+2}{N}}}{ |x-y|^{N-\alpha}}  dxdy \right.  \\
    &\left.-\int_{\mathbb{R}^N} \int_{\mathbb{R}^N}  \frac{(\frac{h_\lambda(e^{Ns/2} u(x))}{2}- \lambda)| u(x)|^{1+\frac{\alpha+2}{N}}h_\lambda(e^{Ns/2}u(y))|u(y)|^{1+\frac{\alpha+2}{N}} }{ |x-y|^{N-\alpha}}  dxdy \right]  \\
  \rightarrow & -\infty, \   \text{as} \ s  \rightarrow +\infty,
\end{align*}
which completes the proof.
\hfill$\Box$

{\bf Remark 2.4.} {\it For any $t\neq 0 ,$ we have
 $f (t)t > \frac{N+\alpha}{N}F(t). $\\
Indeed, define
$$
g(t)=
\begin{cases} \frac{f(t)t- \frac{N+\alpha}{N}F(t)}{|t|^{1+(\alpha+2)/N}}, & t\neq0,\\
0, & t=0.
\end{cases}
$$
By $(H0)-(H1) $ and $(H3),$ we obtain
that $f (t)t > \frac{N+\alpha}{N}F(t) $ for any $t\neq 0. $
}

{\bf Lemma 2.5.} {\it   Let  $(H0)-(H3)$ be  satisfied, then\\
(i)  $F(t) > 0 $  for any $t \neq 0.$\\
(ii) There exists a positive sequence $\{\tau_n^+\}$ and a negative sequence  $\{\tau_n^-\}$
such
that  $|\tau_n^\pm| \rightarrow 0$ and
$$f (\tau_n^\pm)\tau_n^\pm
 > \frac{N+\alpha+2}{N}F(\tau_n^\pm), \  \text{for each } n \geq 1.$$
(iii) There exists a positive sequence $\{\sigma_n^+\}$ and a negative sequence $\{\sigma_n^-\}$
 such that $| \sigma_n^\pm  | \rightarrow +\infty $  and
$$f (\sigma_n^\pm )\sigma_n^\pm > \frac{N+\alpha+2}{N}F(\sigma_n^\pm ),\  \text{for each } n \geq 1.$$
(v) For any $t \neq 0,$
 $$f (t)t > \frac{N+\alpha+2}{N}F(t) > 0 .$$}
\quad{\bf Proof.}
(i) Assume by contradiction that $F(t_0) \leq 0$ for some  $t_0 \neq0. $
By  $(H0)-(H1),$  the function  $F(t)/|t|^{1+(\alpha+2)/N}$
reaches its global minimum at some $\tau \neq  0$ such that  $F(\tau )\leq 0 $  and
$$
\left[\frac{F(t)}{|t|^{1+(\alpha+2)/N}} \right]'_{t=\tau}
=\frac{f (\tau)\tau - \frac{N+\alpha+2}{N}F(\tau )}{|\tau|^{2+(\alpha+2)/N}sgn(\tau)}
= 0.\eqno(2.6)$$
Hence, by Remark 2.4, we have
$$0 < f (\tau )\tau -\frac{N+\alpha}{N}F(\tau ) = \frac{2}{N}F(\tau) \leq 0. \eqno(2.7)$$
 This is a contradiction.

(ii)
First, we  consider the positive case.
Assume by contradiction that there exists $T_s > 0$
small enough such that
$$f (t)t \leq \frac{N+\alpha+2}{N}F(t) \ \text{ for any }  t \in(0, T_s ].$$
By (i) and (2.6), we have
$F(t)/|t|^{1+(\alpha+2)/N}  \geq  F(T_s)/|T_s|^{1+(\alpha+2)/N} >0 $ for all  $t \in(0, T_s ].$
Thus,  $\lim_{t \rightarrow 0} F(t)/|t|^{1+\frac{(\alpha+2)}{N}}> 0$ which contracts $(H1).$
The proof for the negative case is similar.

(iii)
Noting that the two cases are similar, we only show the existence of $\{\sigma_n^-\}$ .
By contradiction we suppose
that there exists $T_l > 0$ such that
 $$f (t)t \leq \frac{N+\alpha+2}{N} F(t) \text{ for any} \  t \leq -T_l.$$
This yields that
$F(t)/|t|^{1+(\alpha+2)/N}   \leq  F(-T_l )/|-T_l|^{1+(\alpha+2)/N}  < +\infty $ for all $t < -T_l ,$
and then $\lim _{t \rightarrow  -\infty} F(t)/|t|^{1+(\alpha+2)/N}   < +\infty $   which is in  contradicts with  $(H1).$
 Consequently, the sequence  $\{\sigma_n^-\}$ exists and we have the desired conclusion.

(v)  First we  claim that  $f (t )t \geq \frac{N+\alpha+2}{N} F(t )$ for any $t \neq 0.$

By contradiction we assume that $f (t_0)t_0 < \frac{N+\alpha+2}{N} F(t_0)$ for some $t_0\neq 0.$
Being the cases $t_0 < 0 $ and $t_0> 0$  similar, we   only study the case that $t_0 < 0.$
In view of (ii) and (iii), there exist $\tau_{min}, \tau_{max} \in \mathbb{R} $ such that $\tau_{min} < t_0 < \tau_{max} < 0$
$$f (t)t <\frac{N+\alpha+2}{N}F(t) \ \text{ for any} \ t \in  (\tau_{min}, \tau_{max}  ),\eqno(2.8)$$
and
$$f (t)t = \frac{N+\alpha+2}{N}F(t) \ \text{  when} \   t =\tau_{min}, \tau_{max}. \eqno(2.9)$$
By (2.6) and (2.8), we have
$$F(\tau_{min})/
|\tau_{min}|^{1+(\alpha+2)/N} <
F(\tau_{max })/
|\tau_{max}|^{1+(\alpha+2)/N}  . \eqno(2.10)$$
It follows from (2.9) and  $(H3).$  that
$$\frac{ F(\tau_{min})}{
|\tau_{min}|^{1+(\alpha+2)/N}}
= \frac{N}{2}\frac{\widetilde{F}(\tau_{min})}{|\tau_{min}|^{1+(\alpha+2)/N}}
> \frac{N}{2}\frac{\widetilde{F}(\tau_{max})}{|\tau_{max}|^{1+(\alpha+2)/N}}
=\frac{ F(\tau_{max})}{
|\tau_{max} |^{1+(\alpha+2)/N} },$$
which contradict (2.10). We have the desired conclusion.

By Claim and (2.6), the function $F(t)/|t|^{1+(\alpha+2)/N}$ is nonincreasing on $(-\infty, 0)$ and nondecreasing
on $(0,\infty).$ Then,  the function $f (t)/|t|^{ (\alpha+2)/N}$ is strictly increasing on $(-\infty, 0)$
and $(0,\infty) $  thanks to  $(H3).$
For every  $t \neq 0$,  we see that
$$\frac{N+\alpha+2}{N}F(t) = \frac{N+\alpha+2}{N}
\int_0
^ t f (s)ds
<  \frac{N+\alpha+2}{N}
f (t)/|t|^{ (\alpha+2)/N} \int_0^t |s|^{ (\alpha+2)/N} ds = f (t)t$$
and we obtain the conclusion.
\hfill$\Box$

{\bf Lemma 2.6.} {\it   Let   $(H0)-(H3)$  be satisfied.
Then, for any $u\in H^1(\mathbb{R}^N ) \backslash \{0\},$ \\
(i ) There exists a unique number $s(u) \in \mathbb{R}$  such that $P(s(u)\star u) = 0.$\\
(ii ) $I (s(u)\star u) > I (s\star u)$  for any $s \neq s(u). $ In particular, $I (s(u)\star u) > 0.$\\
(iii ) The mapping $u \mapsto  s(u)$ is continuous in $u \in  H^1(\mathbb{R}^N ) \backslash \{0\}.$\\
(iv) $s(u(\cdot  + y)) = s(u)$  for any $y  \in \mathbb{R}^N . $
 If $f$ is odd, then   $s(-u) = s(u).$
}

{\bf Proof.}
(i ) Noting  that
$$I (s \star u) = \frac{e^{2s}}{2}\int_{\mathbb{R}^N} | \nabla u|^2dx-\frac{e^{-(N+\alpha )s}}{2}\int_{\mathbb{R}^N} \left(I_\alpha  \ast F(e^{Ns/2}u)\right)F(e^{Ns/2}u)dx,$$
we have that $I (s \star u) $ is of class $C^1$ and
$$\frac{d}{ds}I (s \star u) =e^{2s}\int_{\mathbb{R}^N} | \nabla u|^2dx - \frac{e^{-(N+\alpha )s}}{2}N\int_{\mathbb{R}^N} \left(I_\alpha  \ast F(e^{Ns/2}u)\right)\widetilde{F}(e^{Ns/2}u)dx = P(s \star u).$$
From Lemma 2.3, it follows that $ I (s \star u)$ reaches its global maximum at some $s(u) \in \mathbb{R} $ and then
$P(s(u)\star u) = \frac{d}{ds}I (s \star u) = 0.$
 Since
$\widetilde{F}(t) = g(t)|t|^{1+(\alpha+2)/N}$
due to (2.5),
we conclude that
$$P(s\star u) = e^{2s}\left(\int_{\mathbb{R}^N} | \nabla u|^2dx - \frac{e^{ \frac{-(N+\alpha+2  )s}{2}} }{2}N\int_{\mathbb{R}^N} \left(I_\alpha  \ast { F(e^{Ns/2}u)} \right)g(e^{Ns/2}u) |u|^{  \frac{N+\alpha+2}{N}} dx\right).
$$
Fixing  $t\in \mathbb{R}  \backslash  \{0\},$ by $ (H3)$ and Lemma 2.5 (v),
the functions  $s\mapsto g(e^{Ns/2}t)$ and $s\mapsto F(e^{Ns/2}t)$ are strictly increasing due to the fact that
$F(e^{Ns/2}t)  =\frac{ F(e^{Ns/2}t)}{   |e^{Ns/2}t|^  {\frac{N+\alpha+2}{N} } } |e^{Ns/2}t|^  \frac{N+\alpha+2}{N}   $.
We then  have that  $s(u)$ is unique.

(ii ) This follows from (i).

(iii ) From    Item (i ), the mapping $u \mapsto s(u)$ is well-defined.
 Let $u \in H^1(\mathbb{R}^N ) \backslash  \{0\}$   and
$\{u_n\} \subset  H^1(\mathbb{R}^N )  \backslash  \{0\}$ be any sequence such that $ u_n \rightarrow u $ in $H^1(\mathbb{R}^N ).$
Let $ s_n := s(u_n)$ for any $n \geq 1.$
It suffices  to prove that up to a subsequence $s_n \rightarrow  s(u)$  as  $ n \rightarrow\infty .$

We first show that $\{s_n\}$ is bounded. Recall the continuous coercive function $h_\lambda$ defined by
(2.3).
Assume by contradiction that up to a subsequence $s_n \rightarrow  +\infty,$  by
Fatou's lemma and the fact that $u_n  \rightarrow  u \neq 0$ almost everywhere in $\mathbb{R}^N$, we have
$$   \lim_{n \rightarrow \infty} \int_{ \mathbb{R}^N} (I_\alpha  \ast [h_0(e^{\frac{Ns_n}{2} }u_n)|u_n|^{1+\frac{\alpha+2}{N}}])   h_0(e^{Ns_n/2 }u_n)|u_n|^{1+\frac{\alpha+2}{N}}dx =+\infty. $$
From Item (ii ),  we then obtain\\
\begin{align*}
 0\leq e^{-2s_n} I(s_n\star u_n)
 = &\frac{1}{2}\int_{\mathbb{R}^N} | \nabla u_n|^2dx\\
 &-\frac{1}{2}\int_{\mathbb{R}^N} (I_\alpha  \ast [h_0(e^{\frac{Ns_n}{2} }u_n)|u_n|^{1+\frac{\alpha+2}{N}}])h_0(e^{\frac{Ns_n}{2} }u_n)|u_n|^{1+\frac{\alpha+2}{N}}dx\\
 \rightarrow & -\infty ,
\end{align*}
which is a contradiction.
Thus, the sequence $\{s_n\}$ is bounded from above.

 By Item (ii ), one has
$$I (s_n \star u_n) \geq I (s(u)\star u_n) \  \text{for any}\  n \geq 1.$$
Since $s(u)\star u_n \rightarrow s(u)\star u$ in $H^1(\mathbb{R}^N)$, together with Lemma 2.8, we have
$$I (s(u)\star u_n) = I (s(u)\star u) + o_n(1),$$
and there holds
$$\liminf_{n \rightarrow \infty } I (s_n\star u_n)  \geq  I (s(u)\star u) > 0. \eqno(2.11)$$
Note that  $\{s_n \star u_n\} \subset B_m$ for $m > 0$ large enough. Together with Lemma 2.2 (i ) and the fact,
$$\|  \nabla (s_n \star u_n)\|_2= e^{s_n}\| \nabla u_n\|_2,$$
  we deduce from (2.11) that $\{s_n\}$ is bounded also from below.
  In fact,   assume by contradiction  that $s_n \rightarrow -\infty,$ $\lim_{n \rightarrow \infty} I(s_n \star u_n)=0 $ due to Lemma 2.2 (i) which contracts with (2.11).
Without loss of generality,   assume that
$s_n \rightarrow s^* $ for some $s^*\in \mathbb{R}.$
Since $ u_n \rightarrow  u$ in $H^1(\mathbb{R}^N),$  one then has $ s_n \star u_n  \rightarrow s^* \star u $ in $H^1(\mathbb{R}^N).$
Recalling  that
$P(s_n \star u_n) = 0$ for any $n \geq 1,$  from Lemma 2.8 and Lemma 5.5, it follows that
$P(s^* \star u ) = 0.$
Item (i ) implies  that $ s^* = s(u)$ and we have the desired conclusion Item (iii ).

(iv) For any $y \in \mathbb{R}^N,$ after changing variables in the integrals, we have
$P(s(u)\star u(\cdot + y)) = P(s(u)\star u) = 0$
and thus $s(u(\cdot + y))$ = $s(u)$ via Item (i ).
Suppose that  $f$ is odd, we have
$P(s(u)\star (-u)) = P(-(s(u)\star u)) = P (s(u)\star u)  = 0$
and hence $s(-u) = s(u).$
 \hfill$\Box$

{\bf Lemma 2.7.} {\it  Let  $(H0)-(H3)$  be  satisfied. Then \\
(i ) $ \mathcal{P}_m \neq\emptyset  ,$\\
(ii )  $\inf_{u\in \mathcal{P}_m} \|\nabla u \|_{2 }> 0,$\\
(iii ) $E_m:= \inf_{u\in \mathcal{P}_m}  I (u) > 0,$\\
(iv) $I $ is coercive on $\mathcal{P}_m,$ that is $I (u_n) \rightarrow +\infty$
 for any $\{u_n\} \subset \mathcal{P}_m$ with $ \| u_n\|_{H^1} \rightarrow +\infty.$}

{\bf Proof.}
(i )   Item (i) follows from Lemma 2.6 (i ).

(ii ) Assume by contradiction that there exists $\{u_n\}  \subset  \mathcal{P}_m $ such that $\| \nabla u_n\|_2 \rightarrow  0, $
then  by similar arguments as Lemma 2.2 (i), we have
$$ \int_{\mathbb{R}^N} (I_\alpha \ast F(u_n)) \widetilde{F}(u_n) dx \leq \frac{1}{N}  \| \nabla u_n\|_2^2, $$
and thus
$$0 = P(u_n) \geq  \frac{1}{2}
\int_{\mathbb{R}^N}
|\nabla u_n|^2dx >0$$
for $n$ large enough,
which is a contradiction.
Therefore, $\inf_{u\in  \mathcal{P}_m }
\int_{\mathbb{R}^N}
|\nabla u_n|^ 2dx > 0.$

(iii ) For any $u\in  \mathcal{P}_m , $ by Lemma 2.6 (i ) and (ii ), we have
$I (u) = I (0\star u) \geq  I (s\star u)$ for all $s \in \mathbb{ R}.$
Choose $\delta> 0$ be the number given by Lemma 2.2 (i ) and
$s := \ln( \delta/\|\nabla u \|_2) . $
Since
$\|\nabla (s\star u )\|_2=\delta$, from Lemma 2.2 (i ), we have
$$I (u) \geq  I (s\star u) \geq  \frac{1}{4} \int_{\mathbb{R}^N} |\nabla (s\star u ) |^2 = \frac{1}{4}\delta^2 $$
and we obtain the conclusion  Item (iii ).

(iv) Assume by contradiction  that there exists $\{u_n\}  \subset  \mathcal{P}_m $  such that
$\| u_n \|_{H^1}\rightarrow \infty  $
and $\sup_{ n\geq 1} I (u_n) \leq  c$ for some $c \in  (0,+\infty).$
 For any $n \geq 1,$ define
$s_ n := \ln(\|\nabla u_n \|_2)$
 and $v_n := (-s_n)\star u_n.$
Clearly, $s_n \rightarrow +\infty, $ $\{v_n\}  \subset  S_m $ and
 $\|\nabla v_n \|_2
= 1 $ for any $n \geq 1.$ Set
$$\rho  := \limsup_{  n \rightarrow\infty }\left(\sup_{y\in  \mathbb{R}^N}
\int _ { B(y,1) }|v_n|^2dx\right),$$
and we discuss in  two cases.

 Non-vanishing:   $\rho> 0.$ Up to a subsequence, there exists $\{y_n\}\in \mathbb{R}^N$ and
$w \in  H^1(\mathbb{R}^N )\backslash  \{0\}$ such that
$w_n := v_n(\cdot + y_n) \rightharpoonup w $ in $H^1(\mathbb{R}^N )$ and $w_n \rightarrow  w$ a.e. in $\mathbb{R}^N $.
From the fact  $ s_n  \rightarrow +\infty $
 and Fatous lemma, it follows that
 $$   \lim_{s_n \rightarrow +\infty}\int_{\mathbb{R}^N} (I_\alpha  \ast [ h_0(e^{\frac{Ns_n}{2} }w_n)|w_n|^{1+\frac{\alpha+2}{N}}]) h_0(e^{Ns_n/2 }w_n)|w_n|^{1+\frac{\alpha+2}{N}} =+\infty. $$
  Item (iii )  gives that
$$
\begin{array}{rl} 0&\leq e^{-2s_n} I(s_n\star v_n)\\
&= \frac{1}{2}\int_{\mathbb{R}^N} | \nabla v_n|^2dx
-\frac{1}{2}\int_{\mathbb{R}^N} (I_\alpha  \ast [ h_0(e^{\frac{Ns_n}{2} }v_n)|v_n|^{1+\frac{\alpha+2}{N}}])h_0(e^{\frac{Ns_n}{2} }v_n)|v_n|^{1+\frac{\alpha+2}{N}}dx\\
&\rightarrow -\infty ,
\end{array}
\eqno(2.12)
$$
which is a contradiction.

 Vanishing:  $\rho= 0.$ Applying  Lions Lemma \cite[Lemma I.1]{Lions2}, we deduce that
 $v_n  \rightarrow  0$  in
$L^{2(N+\alpha+2)/(N+\alpha)} (\mathbb{R}^N).$
By Lemma 2.2 (ii ), we thus have
\begin{align*}
  \lim_{n \rightarrow \infty}e^{ -(N+\alpha)s}\int_{\mathbb{R}^N} (I_\alpha  \ast F(e^{Ns/2} v_n))F(e^{Ns/2}v_n)dx =0 \  \text{for any} \  s \in \mathbb{R}.
\end{align*}

Noting that
$P(s_ n \star v_n) = P(u_n) = 0,$  in view of  Lemma 2.6 (i ) and (ii ) we deduce    that, for any $s \in \mathbb{ R},$
\begin{align*}
  c\geq I (u_n) = &I (s_ n \star v_n) \geq I (s  \star v_n)\\
= &\frac{e^{2s}}{2}\int_{\mathbb{R}^N} | \nabla v_n|^2dx-\frac{e^{ -(N+\alpha)s}}{2}\int_{\mathbb{R}^N} (I_\alpha  \ast F(e^{Ns/2} v_n))F(e^{Ns/2} v_n)dx\\
=& \frac{e^{2s}}{2}+ o_n(1).
\end{align*}
which is a contradiction for $s > \ln(2c)/2.$\\
Therefore, $I$ is coercive on $\mathcal{P}_m.$
\hfill$\Box$

 {\bf Lemma 2.8.} {\it   Let $(H0)-(H3)$  hold
 and $ \{u_n\} \subset  H^1(\mathbb{R}^N )$ be bounded  such that
$ u_n  \rightarrow  u$ almost everywhere in  $\mathbb{R}^N.$
  Then
$$  \lim_{ n \rightarrow  \infty} \int_{\mathbb{R}^N } \big( (I_\alpha  \ast F(u_n))F(u_n)-(I_\alpha  \ast F(u_n-u))F(u_n-u)-(I_\alpha  \ast F(u))F(u)  \big)dx =0 .$$
}
\quad{\bf Proof.} Since $ \{u_n\} \subset  H^1(\mathbb{R}^N )$ is bounded, up to a subsequence (still denote $\{u_n\}$  ),
we have $u_n \rightharpoonup u$ in $H^1(\mathbb{R}^N ).$
Let $w_n=u_n-u.$
Hence, we obtain that
$\{w_n\}$ is bounded, $w_n\rightharpoonup 0$ in  $H^1(\mathbb{R}^N ),$ and $w_n \rightarrow 0$ a.e. in $\mathbb{R}^N .$
Using $(H0)-(H3)$  and Sobolev embedding theorem, we know $\{F(w_n)\}$  is bounded in $L^{\frac{2N}{N+\alpha}}(\mathbb{R}^N ).$
Recalling that $F$ is continuous, we have $F(w_n(x))\rightarrow0$ a.e. in $\mathbb{R}^N,$
 and then $F(w_n) \rightharpoonup 0$ in $L^{\frac{2N}{N+\alpha}}(\mathbb{R}^N ).$
  Since the Riesz potential $I_\alpha$ defines a linear continuous map from $L^{\frac{2N}{N+\alpha}}(\mathbb{R}^N ) $  to $L^{\frac{2N}{N-\alpha}}(\mathbb{R}^N )$, we have $I_\alpha  \ast F(w_n) \rightharpoonup 0$ in  $L^{\frac{2N}{N-\alpha}}(\mathbb{R}^N ).$
   Therefore
\begin{align*}
\int_{\mathbb{R}^N } \big(I_\alpha  \ast   F(w_n))F(u)dx
\rightarrow  0,
\ \text{and}\
\int_{\mathbb{R}^N } \big(I_\alpha  \ast   F(u))F(w_n)dx
\rightarrow  0.
\end{align*}
By Hardy-Littlewood-Sobolev inequality, one can see that
\begin{align*}
 &\lim_{ n \rightarrow  \infty} \left|\int_{\mathbb{R}^N } \left((I_\alpha  \ast F(u_n))F(u_n)-(I_\alpha  \ast F(u_n-u))F(u_n-u)-(I_\alpha  \ast F(u))F(u) \right)dx \right|\\
 =& \lim_{ n \rightarrow  \infty} \left|\int_{\mathbb{R}^N } \big((I_\alpha  \ast F(w_n+u))F(w_n+u)-(I_\alpha  \ast F(w_n))F(w_n)-(I_\alpha  \ast F(u))F(u)\big) dx\right|\\
 =&\lim_{ n \rightarrow  \infty}\big| \int_{\mathbb{R}^N } \left((I_\alpha  \ast (F(w_n+u)-F(w_n)-F(u)))F(w_n+u)
\right.\\
& \quad \left. +(I_\alpha  \ast F(w_n))(F(w_n+u)-F(w_n)-F(u))
 +(I_\alpha  \ast F(u))(F(w_n+u)-F(w_n)-F(u))\right.\\
&\quad \left. +(I_\alpha  \ast F(w_n))F(u)+(I_\alpha  \ast F(u))F(w_n)  \right) dx\big| \\
 \leq & C_1\|F(w_n+u)-F(w_n)-F(u) \|_{ \frac{2N}{N+\alpha}} \| F(w_n+u)\|_{ \frac{2N}{N+\alpha}} \\
& \quad +C_2\| F(w_n)  \|_{ \frac{2N}{N+\alpha}} \| F(w_n+u)-F(w_n)-F(u)\|_{ \frac{2N}{N+\alpha}} \\
&\quad  +C_3\| F(u) \|_{ \frac{2N}{N+\alpha}} \| F(w_n+u)-F(w_n)-F(u)\|_{ \frac{2N}{N+\alpha}} .
\end{align*}
It suffices to show that $ \int_{\mathbb{R}^N}|  F(w_n+u)-F(w_n)-F(u)  |^{ \frac{2N}{N+\alpha}}  dx \rightarrow 0.$

By $(H0)-(H3),$   and together with Young's inequality, we have
\begin{align*}
 |F(w_n+u)-F(w_n)| \leq& \int_0^1 |f(w_n+tu)u |dt\\
 \leq  &\int_0^1 \varepsilon| w_n+tu|^{\frac{\alpha +2}{N}}|u| +C_\varepsilon | w_n+tu|^{\frac{\alpha +2}{N-2}}|u|dt\\
 \leq & \int_0^1 [\varepsilon | w_n  |^{\frac{\alpha +2}{N}}|u|+\varepsilon |  u|^{1+\frac{\alpha +2}{N}}
+ C_\varepsilon | w_n |^{\frac{\alpha +2}{N-2}}|u|+C_\varepsilon | u |^{1+\frac{\alpha +2}{N-2}}]dt\\
 \leq & C[ \varepsilon | w_n  |^{\frac{ (N+\alpha +2)}{N }}+\varepsilon | u  |^{\frac{ (N+\alpha +2)}{N }}
 +\varepsilon | w_n  |^{\frac{ (N+\alpha )}{N-2}} +C_\varepsilon | u  |^{\frac{ (N+\alpha )}{N-2}} ].
\end{align*}
Thus,
\begin{align*}
    |  F(w_n+u)-F(w_n)-F(u)  |^{ \frac{2N}{N+\alpha}}
      \leq   C[ \varepsilon | w_n  |^{\frac{2(N+\alpha +2)}{N+\alpha}}+\varepsilon | u  |^{\frac{2(N+\alpha +2)}{N+\alpha}}
 +\varepsilon | w_n  |^{\frac{  2N  }{N-2}} +C_\varepsilon | u  |^{\frac{ 2N }{N-2}} ].
\end{align*}
Let
$$H_n(x):=\max\left\{  |  F(w_n+u)-F(w_n)-F(u)  |^{ \frac{2N}{N+\alpha}}  -  C[ \varepsilon | w_n  |^{\frac{2(N+\alpha +2)}{N+\alpha}}  +\varepsilon | w_n  |^{\frac{  2N  }{N-2}}   ],0 \right\}.$$

In view of Sobolev embedding theorem, we have
$$0\leq H_n(x) \leq C( \varepsilon | u  |^{\frac{2(N+\alpha +2)}{N+\alpha}}+ C_\varepsilon | u  |^{\frac{ 2N }{N-2}})\in L^1(\mathbb{R}^N).$$
By Lebesgue dominated convergence theorem, we have
$$  \int_{\mathbb{R}^N} H_n(x)dx \rightarrow 0,  \ \text{as} \ n \rightarrow \infty . \eqno(2.13)$$
From the definition of $H_n(x), $ we have
$$   |  F(w_n+u)-F(w_n)-F(u)  |^{ \frac{2N}{N+\alpha}}  \leq C[ \varepsilon | w_n  |^{\frac{2(N+\alpha +2)}{N+\alpha}}  +\varepsilon | w_n  |^{\frac{  2N  }{N-2}}   ]+H_n(x), $$
which, together with (2.13), implies that
$$ \int_{\mathbb{R}^N}|  F(w_n+u)-F(w_n)-F(u)  |^{ \frac{2N}{N+\alpha}}  dx \rightarrow 0.$$
And we have the assertion.  \hfill$\Box$

\section{Properties of the function $m  \mapsto  E_m$}\label{s3}

{\bf Lemma 3.1.} {\it Let  $(H0)-(H3)$   be satisfied. Then  the function $m  \mapsto  E_m$ is
continuous at each $m > 0.$}

{\bf Proof.}
For a given $ m > 0 ,$  let  $\{m_k \} \subset \mathbb{R}$  satisfying
  $m_k \rightarrow m $ as $k\rightarrow \infty.$
It suffices to show that
$ \lim_{k\rightarrow \infty} E_{m_k}= E_m.$

Firstly, we  show that
$$ \limsup_{k\rightarrow \infty} E_{m_k}\leq E_m.\eqno(3.1)$$

For any $u \in \mathcal{P}_m,$ let
$$u_k :=\sqrt{\frac{m_k}{m}}u \in S_{m_k} , \  k \in \mathbb{N}^+.$$
Since $u_k \rightarrow u$  in $H^1(\mathbb{R}^N ),$
using Lemma 2.6 (iii ),
 we have  $\lim_{k\rightarrow \infty} s(u_k ) = s(u) = 0 $
 and
this yields that
$s(u_k)\star u_k  \rightarrow  s(u)\star u = u $  in $H^1(\mathbb{R}^N ),$   as $k\rightarrow \infty.$
Consequently, by Lemma 2.8,
$$\limsup_{k\rightarrow \infty}E_{m_k} \leq \limsup_{k\rightarrow \infty}I (s(u_k)\star u_k) = I (u).$$
Noting that $ u \in  \mathcal{P}_m$ is arbitrary, we obtain  (3.1).
Next, we  show that
$$\liminf_{k\rightarrow \infty}E_{m_k} \geq E_m. \eqno(3.2)$$
For each $k \in \mathbb{N}^+,$ there exists $v_k \in  \mathcal{P}_{m_k}$ such that
$$I (v_k ) \leq  E_{m_k}+ \frac{1}{k}. \eqno(3.3)$$

Define
$t_k := (\frac{ m}{m_k})^{1/N}$
and $\widetilde{v}_k := v_k (\cdot/t_k ) \in S_m,$
using Lemma 2.6 (ii ) and (3.3), we have
\begin{subequations}
\begin{align}
E_m \leq I (s(\widetilde{v}_k)\star \widetilde{v}_k)
\leq & I(s(\widetilde{v}_k)\star v_k) + |I(s(\widetilde{v}_k)\star \widetilde{v}_k)-I(s(\widetilde{v}_k)\star v_k)|\nonumber\\
\leq & I(  v_k) + |I(s(\widetilde{v}_k)\star \widetilde{v}_k)-I(s(\widetilde{v}_k)\star v_k)| \tag{3.4}\\
\leq & E_{m_k}+ \frac{1}{k} + |I(s(\widetilde{v}_k)\star \widetilde{v}_k)-I(s(\widetilde{v}_k)\star v_k)|.\nonumber
\end{align}
\end{subequations}
It is easy to check that
\begin{align*}
  &|I(s(\widetilde{v}_k)\star \widetilde{v}_k)-I(s(\widetilde{v}_k)\star v_k)| \\
  =&\left| \frac{1}{2}(t_k^{N-2}-1) \| \nabla [s(\widetilde{v}_k)\star v_k] \|_2^2 -\frac{1}{2}(t_k^{N+\alpha}-1)\int_{\mathbb{R}^N} \left(I_\alpha  \ast F(s(\widetilde{v}_k)\star v_k)\right)F(s(\widetilde{v}_k)\star v_k)dx\right|.
\end{align*}

Claim 1. The sequence $ \{v_k \}$ is bounded in $H^1(\mathbb{R}^N ).$

Applying (3.1) and (3.3) and Lemma 2.7, we can deduce that Claim 1 holds.

Claim 2. The sequence $\{ \widetilde{v}_k \} $  is bounded in $H^1(\mathbb{R}^N ),$
 and there exists $\{y_k\}  \subset\mathbb{R}^N $ and
$v\in H^1(\mathbb{R}^N ) $  such that up to a subsequence $\widetilde{v}_k(\cdot +y_k) \rightarrow  v \neq 0$ almost everywhere in $\mathbb{R}^N  .$

By Claim 1, it is easy to see that $\{ \widetilde{v}_k \} $  is bounded in $H^1(\mathbb{R}^N ).$
Set
$$\rho= \limsup_{k \rightarrow \infty} \left(   \sup_{y\in\mathbb{R}^N   }   \int_{B(y,1)} |\widetilde{v}_k |dx \right).$$

  Assume by contradiction that $\rho =0. $
  Then $ \widetilde{v}_k \rightarrow 0$ in $L^{\frac{2(N+\alpha +2)}{N+\alpha}}(\mathbb{R}^N  ).$
 Since $t_k \rightarrow 1,$ we have
  $$\|  {v}_k \|_{\frac{2(N+\alpha +2)}{N+\alpha}}^{\frac{2(N+\alpha +2)}{N+\alpha}}
  = \|  {\widetilde{v}}_k (t_k \cdot )\|_{\frac{2(N+\alpha +2)}{N+\alpha}}^{\frac{2(N+\alpha +2)}{N+\alpha}}
  =t_k^{-N} \|  {\widetilde{v}}_k  \|_{\frac{2(N+\alpha +2)}{N+\alpha}}^{\frac{2(N+\alpha +2)}{N+\alpha}} \rightarrow  0,$$

Note that  $P(v_k)=0,$ and by Lemma 2.2 (ii), we have
$$  \int_{\mathbb{R}^N } |\nabla v_k |^2dx =\frac{N}{2} \int_{\mathbb{R}^N }  (I_\alpha \ast F(v_k)) \widetilde{F}(v_k)dx\rightarrow 0. $$
Thus,
$$0=P(v_k) \geq \frac{1}{2}\int_{\mathbb{R}^N } |\nabla v_k |^2dx  >0$$
which is a contradiction.

Claim 3. $\lim \sup_{k\rightarrow \infty } s(\widetilde{v}_k) < +\infty.$

Assume by contradiction that,  up to a subsequence
$$s(\widetilde{v}_k) \rightarrow +\infty \   \text{as } \ k\rightarrow +\infty. \eqno(3.5)$$
  Claim 2 implies that up to a subsequence
$$\widetilde{v}_k (\cdot +y_k) \rightarrow  v\neq0 \ \text{     a.e. in  } \ \mathbb{R}^N \eqno(3.6)$$
Using Lemma 2.6 (iv) and (3.5), we have
$$s(\widetilde{v}_k (\cdot+y_k)) = s(\widetilde{v}_k )\rightarrow +\infty, \eqno(3.7)$$
and Lemma 2.6 (ii ) implies  that
$$I (s(\widetilde{v}_k (\cdot+y_k))\star \widetilde{v}_k (\cdot+y_k))  \geq 0. \eqno(3.8)$$
Now, in view of (3.6), (3.7) and (3.8), this is a contradiction by using the arguments as (2.12),
and we have the desired conclusion  Claim 3.
Thus, by similar arguments as Lemma 2.7 (i), we have
$$\limsup_{k\rightarrow \infty }  \|s(\widetilde{v}_k (\cdot+y_k)) \star v_k \|_{H^1}   < +\infty.$$
Taking into account $( f _1)$-$( f_3)$ and (3.4), we have the assertion (3.2).
\hfill$\Box$

{\bf Lemma 3.2. } {\it   Let  $(H0)-(H3)$  be satisfied.  Then   the function $m  \mapsto  E_m$ is
nonincreasing on $(0,\infty).$}

{\bf Proof.}
 It   suffices to show that   for every $\varepsilon>0$ and  $m > m'>0$ we have
 $$ E_m \leq E_{m'}+\varepsilon .$$
Take $u\in \mathcal{P}_{m'}$ such that
for any $ \varepsilon>0,$
$$I(u) \leq E_{m'}+\frac{ \varepsilon}{2},
 \eqno(3.9)$$
and define $\chi  \in C_0^\infty( \mathbb{R}^N)$  by
$$\chi (x) =\begin{cases}  1, &|x| \leq 1,\\
 \in [0, 1], &|x| \in (1, 2),\\
 0, &|x| \geq 2.
\end{cases} $$
For every $\delta > 0,$ set $u_{\delta}(x) = u(x) \cdot \chi(\delta x) \in H^1(\mathbb{R}^N ) \backslash \{0\}.$
Note that $u_\delta \rightarrow  u $ in $H^1(\mathbb{R}^N )$
 as $\delta \rightarrow 0^+.$
 By Lemma 2.6 (iii ),
we obtain $\lim_{\delta \rightarrow 0^+} s(u_\delta) = s(u) = 0 $
 and thus
$s(u_\delta)\star u_\delta \rightarrow s(u)\star u = u $ in $H^1(\mathbb{R}^N )$ as $\delta \rightarrow 0^+.$
Consequently, fixing a   $\delta > 0$ small enough, we get
$$I (s(u_\delta)\star u_\delta) \leq I (u) +\frac{\varepsilon}{4}.
\eqno(3.10)$$
Now choose $v \in C_0^\infty(\mathbb{R}^N )$ satisfying supp$(v) \subset  B(0, 1 + 4/\delta) \backslash B(0, 4/\delta)$ and set
$$\widetilde{v}=\frac{m - \|u_\delta\|_2^2 }{ \|v\|_2^2 }v.$$
For any $ \lambda\leq 0, $ let $w_\lambda = u_\delta + \lambda \star \widetilde{v}.$
Choosing the suitable constant $ \lambda,$ we have
supp$(u_\delta) \cap $ supp$(\lambda \star \widetilde{v} ) = \emptyset,$
and thus $w_\lambda\in S_m. $

Claim: $s(w_\lambda)$ is bounded from above when $\lambda \rightarrow -\infty.$

 Assume by contradiction that $s(w_\lambda)\rightarrow +\infty $ as $\lambda \rightarrow -\infty.$
Since $I (s(w_\lambda) \star w_\lambda)\geq  0 $ by Lemma 2.6 (ii )
and that $w_\lambda \rightarrow u_\delta \neq 0$  almost everywhere
in $\mathbb{R}^N $ as  $ \lambda \rightarrow -\infty,$
 we have a contradiction by using similar arguments as (2.12).

Noting that
$$s(w_\lambda ) +\lambda \rightarrow -\infty , \ \text{as} \ \lambda \rightarrow -\infty,$$
we have
$$  \| \nabla (s(w_\lambda ) +\lambda ) \star \widetilde{v}\|_2  \rightarrow 0 .$$

From Lemma 2.2 (ii ), it  follows that
$$  I( ( s(w_\lambda ) +\lambda ) \star \widetilde{v} ) \leq \frac{ \varepsilon}{4} \  \text{for } \ \lambda < 0\ \text{ small enough.}
\eqno (3.11)$$
By virtue of  Lemma 2.6 (ii ), (3.9)-(3.11), we have
\begin{align*}
E_m \leq  I (s(w_\lambda)\star  w_\lambda)
= &I (s(w_\lambda)\star u_\delta)+I (s(w_\lambda)\star  ( \lambda  \star \widetilde{v} ))\\
\leq &I (s(u_\delta)\star u_\delta)+I ((s(w_\lambda)+\lambda)\star    \widetilde{v}  )\\
\leq & I (u) + \frac{\varepsilon}{2}\\
\leq & E_{m'}+ \varepsilon.
\end{align*}
Hence, we obtain the conclusion.\hfill$\Box$

 {\bf Lemma 3.3.} {\it   Let  $(H0)-(H3)$  be satisfied.
 Let $u \in S_m$
and $\mu \in  \mathbb{R}$ such that
$$-\triangle u + \mu u = (I_\alpha \ast F(u))f(u),$$
and  $I (u) = E_m. $  Then $E_m > E_{m'}$ for any $ m' > m$   close enough to $m$  if $\mu > 0$ and for each
$m'< m $ near enough to $ m$  if  $\mu < 0.$}

{\bf Proof.}
For every $t > 0$ and $s \in \mathbb{R},$ we define $u_{t,s} := s\star (tu) \in S_{mt^2}.$
Noting that
$$ \alpha(t,s):= I(u_{t,s})=
  \frac{t^2 e^{2s}}{2}\int_{\mathbb{R}^N} | \nabla u|^2dx-\frac{e^{-(N+\alpha)s}}{2}\int_{\mathbb{R}^N} (I_\alpha  \ast F(e^{Ns/2}tu))F(e^{Ns/2}tu)dx,$$
we have
  $$ \frac{\partial}{\partial t}\alpha(t,s):=
   t e^{2s} \int_{\mathbb{R}^N} | \nabla u|^2dx- e^{-(N+\alpha)s} \int_{\mathbb{R}^N} (I_\alpha  \ast F(e^{Ns/2}tu))f(e^{Ns/2}tu)e^{Ns/2}udx=t^{-1} I'(u_{t,s}) u_{t,s} .$$
For  $\mu> 0,$
$$I'(u)u = - \mu   \|u\|_2^2= -\mu m < 0.$$
Noting that  $u_{t,s} \rightarrow  u$ in $H^1(\mathbb{R}^N ) $
as $(t, s) \rightarrow  (1, 0),$ fixing  $\delta  > 0$ small enough, we have
$$\frac{\partial}{\partial t}\alpha(t,s) <0  \ \text{for any}  \  (t, s) \in (1, 1 + \delta] \times [-\delta, \delta].$$
By the mean value theorem, we have that
$$\alpha(t, s) =\alpha(1, s) + (t - 1) \frac{\partial}{\partial t}\alpha(\theta,s) < \alpha(1, s), \eqno(3.12)$$
where $1 < \theta < t \leq 1 + \delta $  and $|s|\leq \delta.$
By Lemma 2.6 (iii ),   $s(tu) \rightarrow s(u) = 0$ as $t\rightarrow1^+. $
 For any $m'> m$ close enough to $m, $ set
$t :=
\sqrt{\frac{m'}{m}} \in (1, 1 + \delta ]$ and $s := s(tu) \in [-\delta,\delta].$
Thus, from (3.12) and Lemma 2.6 (ii ), it follows that
$$E_{m'}\leq  \alpha(t, s(tu)) < \alpha(1, s(tu)) = I (s(tu)\star u) \leq  I (u) = E_m.$$
The case  $\mu < 0$ can be proved similarly.
 \hfill$\Box$

{\bf Lemma 3.4.} {\it  Let  $(H0)-(H3)$  be satisfied. Let $u \in S_m$
and $\mu \in  \mathbb{R}$ such that
$$-\triangle u + \mu u = (I_\alpha \ast F(u))f(u),$$
and  $I (u) = E_m,$ then   $\mu \geq  0$ .
If in addition $\mu > 0,$  then  $E_m > E_{m'}$ for any $ m' > m.$ }

 {\bf Proof.} From Lemmas 3.2 and 3.3, we directly obtain the conclusion.\hfill$\Box$

{\bf Lemma 3.5.} {\it   Let $(H0)-(H3)$   be satisfied.
 Then $E_m\rightarrow +\infty$ as  $m \rightarrow 0^+.$ }

 {\bf Proof.}
 It suffices to show that for every sequence $\{u_n\}  \subset  H^1(\mathbb{R}^N ) \backslash \{0\} $
 satisfying
$P(u_n) = 0$ and $\lim_{n \rightarrow \infty } \|u_n \|_2= 0,$
it must be  $I (u_n)\rightarrow +\infty $ as $n\rightarrow \infty.$

 Set
$s_n := \ln( \| \nabla u_n \|_2) $ and $v_n := (-s_n)\star u_n.$
Clearly, $\| \nabla v_n \|_2= 1 $  and
$\|   v_n \|_2 \rightarrow  0. $
By Lions Lemma,  $v_n \rightarrow  0 $ in $L^{\frac{2(N+\alpha+2)}{N+\alpha}}  (\mathbb{R}^N).$
 Applying  Lemma 2.2 (ii ), we have
$$\lim_{n \rightarrow \infty}e^{-(N+\alpha)s} \int_{ \mathbb{R}^N}(I_\alpha \ast  F(e^{Ns/2}v_n))F(e^{Ns/2}v_n)dx =0 $$
for any  $s \in \mathbb{R}.$
Observing  that $P(s_n \star v_n) = P(u_n) = 0,$ with the aid  of Lemma 2.6 (i ) and (ii ), we derive
$$I (u_n) = I (s_n \star v_n) \geq I (s \star v_n)
= \frac{1}{2}e^{2s}-\frac{1}{2} e^{-(N+\alpha)s} \int_{ \mathbb{R}^N}(I_\alpha \ast  F(e^{Ns/2}v_n))F(e^{Ns/2}v_n)dx =\frac{1}{2}e^{2s}+o_n(1).$$
Since $ s \in \mathbb{R}$ is arbitrary, we have that $I (u_n)\rightarrow +\infty.$
 \hfill$\Box$

{\bf Lemma 3.6.} {\it  Let $(H0)-(H4)$  be satisfied.
 Then $E_m\rightarrow 0$ as  $m \rightarrow +\infty.$ }

{\bf Proof.}
Take $u \in  S_1 \cap L^\infty (\mathbb{R}^N)$ and let $u_m =\sqrt{m}u \in S_m$ for any $m > 1.$
 By Lemma 2.6 (i ),
we can find a  unique $s(m) \in \mathbb{R} $ such that $s(m)\star u_m \in \mathcal{P}_m.$
Noting that $F$ is nonnegative by Lemma 2.5,
we  have
$$0 < E_m \leq I (s(m)\star u_m )\leq\frac{1}{2} me^{2s(m)} \int_{\mathbb{R}^N} | \nabla u|^2dx.  \eqno(3.13)$$
 Now, it  suffices to show that
$$\lim_{m \rightarrow +\infty }\sqrt{ m}e^{s(m)}=0.
$$
 By (2.6),  we have
 $\frac{F(t)}{ |t|^{\frac{N+\alpha+2}{N}} }$ is  strictly increasing on $(0,\infty)$ and is strictly decreasing on $(-\infty,0).$
 Combing with $(H1)$, we then have
 $ h_0(t)> 0,$ for $t \neq 0.$

Remembering the functions $g$ and $h_\lambda$ defined by (2.5) and (2.3), by $P(s(m)\star u_m ) = 0, $  we obtain
$$     \int_{\mathbb{R}^N} | \nabla u|^2dx
=   \frac{N}{2} m^{\frac{\alpha+2 }{ N}} \int_{\mathbb{R}^N} (I_\alpha  \ast [ h_0( \sqrt{ m}e^{Ns(m)/2}u) |u|^{\frac{N+\alpha+2}{N}} ] )g(\sqrt{ m}e^{Ns(m)/2}u) |u|^{\frac{N+\alpha+2}{N}}dx. $$
 This gives
$$ \lim_{m  \rightarrow +\infty} \sqrt{ m}e^{Ns(m)/2} =0.$$
Lemma 2.6 and $(H4) $ yield that there exists  $\delta> 0 $ small enough such that
$$
\widetilde{F}(t) \geq  \frac{2}{N}F(t) \geq  \varepsilon^{-1} |t|^{ \frac{ N+\alpha}{N-2}}  \  \text{for any} \ |t| \leq \delta.\eqno(3.14)$$
 From $P(s(m)\star u_m) = 0$ and (3.14), it follows that
 \begin{align*}
  \int_{\mathbb{R}^N} | \nabla u|^2dx
=&m^{-1}\frac{e^{-(N+2+\alpha)s(m)}}{2}N\int_{\mathbb{R}^N} (I_\alpha  \ast F(\sqrt{m}e^{Ns(m)/2}u))\widetilde{F}(\sqrt{m}e^{Ns(m)/2}u)dx\\
\geq & \frac{2}{N}m^{-1}\frac{e^{-(N+2+\alpha)s(m)}}{2} |\sqrt{m}e^{Ns(m)/2}|^{ \frac{2(N+\alpha)}{N-2}}  N\int_{\mathbb{R}^N} (I_\alpha  \ast (\varepsilon^{-1} |u|^{ \frac{ N+\alpha}{N-2}})) (\varepsilon^{-1} |u|^{ \frac{ N+\alpha}{N-2}})dx\\
\geq & \varepsilon^{-2}  {m}^{ \frac{ 2+\alpha}{N-2} } e^{\frac{ 2(2+\alpha)s(m)}{N-2}}\int_{\mathbb{R}^N} (I_\alpha  \ast   |u|^{ \frac{ N+\alpha}{N-2}} ) (|u|^{ \frac{ N+\alpha}{N-2}})dx\\
\geq &  \varepsilon^{-2} ( {\sqrt{m}e^{s(m)}}) ^{ \frac{2( 2+\alpha)}{N-2} } \int_{\mathbb{R}^N} (I_\alpha  \ast   |u|^{ \frac{ N+\alpha}{N-2}} ) (|u|^{ \frac{ N+\alpha}{N-2}})dx
\ \   \text{for large enough} \ m,
 \end{align*}
which implies that   $ \lim_{m \rightarrow \infty } {\sqrt{m}e^{s(m)}}=0$
and then the assertion  holds due to (3.13).
 \hfill$\Box$

\section{Ground states}\label{s3}

For given $m > 0,$ we introduce the   functional
$$\Psi(u):=I(s(u)\star u) =\frac{e^{2s(u)}}{2}  \int_{ \mathbb{R}^N }|\nabla u|^2dx -\frac{e^{-(N+\alpha)s(u) }}{2}
\int_{ \mathbb{R}^N }( I_\alpha \ast F(e^{Ns(u)/2}u ))F(e^{Ns(u)/2}u )dx.$$
\quad{\bf Lemma 4.1.} {\it  The functional $  \Psi : H^1(\mathbb{R}^N ) \backslash  \{0\} \rightarrow  \mathbb{R} $
is of class $C^1$ and
\begin{align*}
d\Psi(u)[\varphi] = &e^{2s(u)} \int_{ \mathbb{R}^N }\nabla u \nabla \varphi dx
-e^{-(N+\alpha)s(u)}\int_{ \mathbb{R}^N } ( I_\alpha \ast F(e^{Ns(u)/2}u ))f (e^{Ns(u)/2}u)e^{Ns(u)/2} \varphi dx\\
= &d I(s(u)\star u)[s(u)\star \varphi ]
\end{align*}
for any  $u \in H^1(\mathbb{R}^N ) \backslash  \{0\}$ and  $\varphi \in H^1(\mathbb{R}^N ) .$
}

{\bf Proof.}
The proof for the Schrodinger equation is given in \cite{Jeanjean4}. Only some adjustments are needed, and we omit it.
\hfill$\Box$

We  consider the constrained
 functional $J =\Psi|_{S_m} : S_m  \rightarrow \mathbb{R}$  which is of class $ C^1$ and
$$d J(u)[\phi] = d\Psi(u)[\phi] = d I(s(u)\star u)[s(u)\star \phi]$$
for any $u \in S_m$ and $ \phi \in  T_u S_m.$

{\bf Lemma 4.2.} {\it There exists a Palais-Smale sequence $\{u_n\} \subset \mathcal{P}_m$ for the constrained functional
$I|_{S_m} $ at the level  $E_m.$
When $f$ is odd, we have in addition $\| u_n^-\|_2\rightarrow 0,$ where $v^-$ stands
for the negative part of $v.$}

{\bf Proof.} By Lemma 4.4, we choose $\mathcal{G}$ the class of all singletons in $S_m.$
Since $E_m>0$ due to Lemma 2.7, it suffices to show that $E_m= E_{m, \mathcal{G}}.$
Note that
$$E_{m, \mathcal{G}} =  \inf_{A \in \mathcal{G}} \max_{u \in A} J(u) =\inf_{ u \in S_m} I(\rho(u) \star u).$$
For any $u\in S_m,$ $\rho(u) \star u  \in \mathcal{P}_m,$
and $I(\rho(u) \star u  ) \geq E_m,$
and then $ E_{m, \mathcal{G}} \geq E_m. $
If $u \in \mathcal{P}_m,$ we obtain $\rho (u)=0$ and $ I(u) \geq E_{m, \mathcal{G}} ,$
thus  $E_m \geq E_{m, \mathcal{G}} .$
\hfill$\Box$

{\bf Definition 4.3. (\cite[ Definition 3.1]{Ghoussoub}).} {\it Let $B$ be a closed subset of a metric space $X.$
 We say
that a class $\mathcal{G}$ of compact subsets of $X$ is a homotopy stable family with closed boundary $ B$
provided\\
(i ) every set in $\mathcal{G}$ contains $B,$\\
(ii ) for any set $A \in \mathcal{G} $ and any homotopy $\eta \in C([0, 1] \times X, X)$ that satisfies $\eta(t, u) = u$
for all $(t, u) \in (\{0\} \times X) \cup ([0, 1] \times B), $
one has $ \eta(\{1\} \times A) \in \mathcal{G}.$ }\\
We remark that the case $B =\emptyset $ is admissible.

{\bf Lemma 4.4.}  {\it  Let $\mathcal{G}$ be a homotopy stable family of compact subsets of $S_m $  (with  $B = \emptyset$) and
set
$$E_{m,\mathcal{G}} := \inf_{A\in \mathcal{G}} \max_{u\in A}J (u).$$
If $E_{m,\mathcal{G}}  > 0,$  then there exists a Palais-Smale sequence $\{u_n\} \subset  P_m$ for the constrained
functional $I|_{S_m}$ at the level $E_{m,\mathcal{G}} .$
 In the particular case when $f$ is odd and $\mathcal{G}$ is the class of
all singletons included in $S_m,$ we have in addition that $\|u_n^-\|_2 \rightarrow 0.$}

{\bf Proof.}
Let $\{A_n\} \subset \mathcal{G} $ be an arbitrary minimizing sequence of $E_{m,\mathcal{G}}.$
We define the map
$$\eta : [0, 1] \times S_m \rightarrow S_m, \ \eta(t, u) = (ts(u))\star u,$$
which is continuous  and well defined due to Lemma 2.6(iii ).
Since
  $\eta(t, u)= u $ for all $(t, u) \in \{0\} \times S_m,$
by the definition of $\mathcal{G},$ one has
$$D_n := \eta(1, A_n) = \{s(u)\star u | u \in A_n\} \in \mathcal{G}. \eqno(4.1)$$
In particular, $D_n \subset  \mathcal{P}_m $ for every $n \in \mathbb{N}^+$.
Since $J (s\star u) = J (u)$ for any $s \in \mathbb{ R}$ and any
$u \in S_m,$ we have
 $$\max_{D_n} J = \max_{A_n} J \rightarrow  E_{m,\mathcal{G}} $$
 and thus $\{D_n\} \subset \mathcal{G}$ is another minimizing
sequence of $E_{m,\mathcal{G}} .$
Now, by the minimax principle \cite[Theorem 3.2]{Ghoussoub} , we can find a Palais-
Smale sequence $\{v_n\} \subset S_m$ for $J$ at the level $E_{m,\mathcal{G}} $ such that dist$_{H^1(\mathbb{R}^N)} (v_n, D_n) \rightarrow 0$ as
$n\rightarrow \infty .$
Define
$$s_n := s(v_n)\ \text{and }  u_n := s_ n \star v_n = s(v_n)\star v_n.$$

Claim. There exists $C > 0 $ such that $ e^{-2s_n}\leq C$ for every $n.$

Noting that
$$e^{-2s_n}=\frac{\| \nabla v_n\|_2^2}{ \| \nabla u_n\|_2^2}.$$
Observing that $\{u_n\} \subset  \mathcal{P}_m,$ from Lemma 2.7 (ii ), it follows that $\{ \| \nabla u_n\|_2^2\} $ is bounded from below by a positive constant.
Since $D_n\subset \mathcal{P}_m $ for every $n,$ we have
$$\max_{D_n}I = \max_{D_n}J \rightarrow E_{m,\mathcal{G}}$$
and clearly, $\{D_n\}$ is uniformly bounded in $H^1(\mathbb{R}^N )$ by Lemma 2.7 (iv) .
Since dist$_{H^1(\mathbb{R}^N)} (v_n, D_n) \rightarrow 0,$   we have
 $\| \nabla v_n\|_2^2 < \infty .$
 Consequently, we obtain the conclusion.

 Now, using  $\{u_n\}\subset \mathcal{P}_m,$  we have
$I (u_n) = J (u_n) = J (v_n)  \rightarrow E_{m,\mathcal{G}}.$
It suffices  to show that $\{u_n\}$ is a Palais-Smale sequence for $I$ on $S_m.$
 For every
$\psi \in T_{u_n} S_m,$ we have
$$  \int_{\mathbb{R}^N } v_n[(-s_n) \star \psi]dx = \int_{\mathbb{R}^N } (s_n \star v_n)\psi dx= \int_{\mathbb{R}^N }u_n \psi dx=0,$$
and it must be $ (-s_n) \star \psi  \in T_{v_n} S_m.$
 Also,  by the Claim, we have  $\| (-s_n) \star \psi \|_{H^1} \leq  \max\{\sqrt{C},1\} \| \psi\|_{H^1}$
Denoting by  $\| \cdot\|_{u,*} $  the dual norm of $(T_u S_m)^*$
 and in view of  the definition of $J|_{S_m}$, we deduce that
\begin{align*}
\|dI(u_n)\|_{u_n,*}
=&\sup_{\psi \in T_{u_n} S_m, \| \psi\|_{H^1} \leq 1 } | dI(u_n) \psi|\\
=&\sup_{\psi \in T_{u_n} S_m, \| \psi\|_{H^1} \leq 1 } | dI(s_n \star v_n)[s_n \star( (-s_n )\star \psi )]|\\
=&\sup_{\psi \in T_{u_n} S_m, \| \psi\|_{H^1} \leq 1 } | dJ(v_n) [(-s_n )\star\psi]|\\
\leq &\| dJ(v_n)\|_{v_n,*} \sup_{\psi \in T_{u_n} S_m, \| \psi\|_{H^1} \leq 1 } \| (-s_n )\star \psi\|_{H^1} \\
\leq & \max\{\sqrt{C},1\}\| dJ(v_n)\|_{v_n,*}.
\end{align*}
Since $\{v_n\}  \subset  S_m$ is a Palais-Smale sequence of $J ,$  we have
$\| dI(u_n)\|_{u_n,*}\rightarrow 0.$

Note that the class of all singletons included in $S_m$ is a homotopy stable family
of compact subsets of $S_m$ (with $B = \emptyset$).
When $f$ is odd, by  the above arguments,  we can
take a minimizing sequence $\{A_n\} \subset \mathcal{G}$
 which consists of nonnegative functions (rather
than an arbitrary one). Thus the sequence $\{D_n\}$ defined in (4.1) inherits this property.
Since
 dist$_{H^1(\mathbb{R}^N)} (v_n, D_n) \rightarrow 0,$ we obtain a Palais-Smale sequence $\{u_n\} \subset \mathcal{P}_m $
 for $I|_{S_m} $ at the level
$E_{m,\mathcal{G}}$ satisfying the additional property
$$  \|u_n^-\|_2^2=\|s(v_n)\star v_n^-\|_2^2 =\|v_n^-\|_2^2 \rightarrow 0.$$
The proof of this lemma is completed.
\hfill$\Box$

{\bf Lemma 4.5.} {\it  Let $\{u_n\} \subset  S_m$ be any bounded Palais-Smale sequence for the constrained
functional $I|_{S_m}$ at the level $E_m > 0 $ satisfying $P(u_n) \rightarrow 0.$
Then there exists $u \in S_m$ and $\mu > 0$ such that,
 up to the extraction of a subsequence and up
to translations in $\mathbb{R}^N,$ $ u_n \rightarrow u $ strongly in $H^1(\mathbb{R}^N)$
 and $ -\triangle u + \mu u = (I_\alpha \ast F(u))f (u).$}

{\bf Proof.}  Since  $\{u_n\} \subset  S_m$ is bounded in $H^1(\mathbb{R}^N)$, up to a subsequence, one may  assume
that $ \lim_{n\rightarrow \infty }  \| \nabla u_n \|_2,$
$ \lim_{n\rightarrow \infty }    \int_{\mathbb{R}^N }(I_\alpha \ast F(u_n ))F(u_n )dx$
and $\lim_{n\rightarrow \infty }    \int_{\mathbb{R}^N }(I_\alpha \ast F(u_n ))f(u_n )u_ndx$  exist.
Recalling that $\| d I(u_n)\|_{u_n,*} \rightarrow 0 $, by \cite[Lemma 3]{Berestycki2}, we have
$$- \Delta u_n + \mu_n u_n -(I_\alpha \ast F(u)) f (u_n) \rightarrow 0  \   \   \text{in} \   (H^1(\mathbb{R}^N ))^*,$$
where
$$\mu_n := \frac{1}{m}\left(  \int_{\mathbb{R}^N }(I_\alpha \ast F(u_n ))f(u_n )u_ndx -\| \nabla u_n \|_2^2 \right).$$
From $\mu_n \rightarrow  \mu $ for some $\mu \in \mathbb{R},$ it follows that
$$- \triangle u_ n(\cdot+ y_n) + \mu u_n(\cdot + y_n) - (I_\alpha \ast F(u_n (\cdot + y_n))) f (u_n(\cdot + y_n)) \rightarrow 0  \  \    \text{in} \    (H^1(\mathbb{R}^N ))^*
\eqno(4.2) $$
for any $\{y_n\} \subset  \mathbb{R}^N  .$

Claim: $\{u_n\} $ is non-vanishing.

 Assume by contradiction that $\{u_n\}$
is vanishing, then $u_n \rightarrow  0 $ in $L^{2(N+\alpha+2)/(N+\alpha)} (\mathbb{R}^N  )$
by Lions Lemma \cite[Lemma I.1]{Lions2}. Noting $P(u_n) \rightarrow 0 $ and Lemma 2.2 (ii ), we obtain
 $\int_{\mathbb{R}^N }(I_\alpha \ast F(u_n ))\widetilde{F}(u_n )dx \rightarrow 0 $  and then
$$  \| \nabla u_n \|_2^2= P(u_n)+\frac{N}{2}\int_{\mathbb{R}^N }(I_\alpha \ast F(u_n ))\widetilde{F}(u_n )dx\rightarrow 0.$$
Hence,
$$E_m = \lim_{n \rightarrow  \infty  }I (u_n) = \lim_{n \rightarrow  \infty  }\left(\frac{1}{2} \| \nabla u_n \|_2^2-\frac{1}{2}\int_{\mathbb{R}^N }(I_\alpha \ast F(u_n ))F(u_n )dx\right)= 0 ,$$
 which is in  contradiction with   $ E_m > 0,$ and so  we obtain the claim.

 Since the sequence $\{u_n\}$ is non-vanishing, up to a subsequence, there exists $\{y_n^1\} \subset  \mathbb{R}^N$
  and $w^1 \in B_m \backslash  \{0\}$ such
that  $u_n(\cdot + y_n^1) \rightharpoonup w^1$ in  $H^1(\mathbb{R}^N ),$
$ u_n(\cdot + y_n^1) \rightarrow  w^1$ in $L_{loc} ^p(\mathbb{R}^N )$ for any $p \in [1, 2^*), $
and $ u_n(\cdot + y_n^1) \rightarrow w^1$  almost everywhere in $\mathbb{R}^N.$
Note that  $|f(u_n(\cdot + y_n^1))|^{ \frac{2N}{N+\alpha}}$ is bounded in $L^{\frac{N+\alpha}{2+\alpha} }(\mathbb{R}^N),$
and  $|f(u_n(\cdot + y_n^1)) - f(w^1) |^{ \frac{2N}{N+\alpha}}\rightharpoonup 0$ in $L^{\frac{N+\alpha}{2+\alpha} }(\mathbb{R}^N).$
Since $|f(u_n(\cdot + y_n^1)) - f(w^1) |^{ \frac{2N}{N+\alpha}}\rightharpoonup 0$ in $L^{\frac{N+\alpha}{2+\alpha} }(\mathbb{R}^N),$ and $I_\alpha \ast F(u_n (\cdot + y_n^1)) \rightharpoonup I_\alpha \ast F(w^1)$ in $ L^ \frac{2N}{N-\alpha}(\mathbb{R}^N),$ we see that
\begin{align*}
&\lim _{n\rightarrow \infty   } \left|\int_{\mathbb{R}^N }[(I_\alpha \ast F(u_n (\cdot + y_n^1))) f(u_n(\cdot + y_n^1) )-(I_\alpha \ast F(w^1)f(w^1) ] \varphi  dx \right|\\
&\leq
\lim _{n\rightarrow \infty   }\int_{\mathbb{R}^N }\left|[(I_\alpha \ast F(u_n (\cdot + y_n^1)))[ f(u_n(\cdot + y_n^1) )-f(w^1)]\varphi \right|dx\\
& +\lim _{n\rightarrow \infty   }\int_{\mathbb{R}^N }|[(I_\alpha \ast F(u_n (\cdot + y_n^1))-(I_\alpha \ast F(w^1)]f(w^1)   \varphi | dx \\
&= 0.
\end{align*}
 We observe from  (4.2) that
$$- \triangle w^1 + \mu w^1 = (I_\alpha \ast F(w^1)) f (w^1). \eqno(4.3)$$
By the Nehari and Pohozaev identities corresponding to (4.3) and Theorem 3 in \cite{Moroz}, we have $P(w^1) = 0.$
Take $v_n^1= u_n-w^1(\cdot -y_n^1)$  for every $ n \in \mathbb{N}^+.$
As a consequence, $v_n^1(\cdot +y_n^1)
\rightharpoonup 0$ in $H^1(\mathbb{R}^N )$ and
$$m = \lim_{n \rightarrow \infty} \| v_n^1(\cdot +y_n^1) + w^1 \|_2^2= \|w^1\|_2^2+ \lim_{n \rightarrow \infty} \| v_n^1\|_2^2
. \eqno(4.4)$$
 By Lemma 2.8, we  have
\begin{align*}
  & \lim_{n \rightarrow \infty}  \int_{\mathbb{R}^N} \left(I_\alpha \ast F(u_n (\cdot + y_n^1))\right) F(u_n (\cdot + y_n^1))dx\\
= &\int_{\mathbb{R}^N} \left(I_\alpha \ast F(w^1)\right) F(w^1)dx
+\lim_{n \rightarrow \infty}  \int_{\mathbb{R}^N} \left(I_\alpha \ast F(v_n^1 (\cdot + y_n^1))\right) F(v_n^1 (\cdot + y_n^1))dx.
\end{align*}
Thus,
$$E_m
= \lim_{n \rightarrow \infty} I(u_n)
 =\lim_{n \rightarrow \infty} I(u_n(\cdot + y_n^1))
= I(w^1) + \lim_{n \rightarrow \infty}I (v_n^1 (\cdot + y_n^1))
= I(w^1) +  \lim_{n \rightarrow \infty}I (v_n^1 ).
\eqno(4.5)$$
Claim: $\lim_{n \rightarrow \infty}I (v_n^1 )\geq 0. $

 Assume by contradiction that
$\lim_{n \rightarrow \infty}I (v_n^1 )< 0. $
Suppose that $\{v_n^1\} $  is vanishing. Then $v_n^1 \rightarrow 0$ in $L^{  \frac{2(N+\alpha+2)}{N+\alpha}}( \mathbb{R}^N),$
and
$\lim_{n \rightarrow \infty}  \int_{\mathbb{R}^N} \left(I_\alpha \ast F(v_n^1\right) F(v_n^1 )dx =0,$
thus,
$$ \lim_{n \rightarrow \infty} I(v_n^1)=  \lim_{n \rightarrow \infty}\frac{1}{2} \left[ \int_{\mathbb{R}^N } |\nabla v_n^1|^2 dx-    \int_{\mathbb{R}^N} \left(I_\alpha \ast F(v_n^1\right) F(v_n^1 )dx \right] \geq 0, $$
which is a contradiction.
Hence, $\{v_n^1\} $  is non-vanishing and, up to a subsequence, there exists a
sequence $\{y_n^2\}\subset \mathbb{R}^N$ such that
$$\lim_{n \rightarrow \infty} \int_{B(y_n^2,1) } |v_n^1|^2dx> 0.$$
Then there holds that $|y_n^2-y_n^1|
\rightarrow\infty   $  (since $v_n^1 (\cdot + y_n^1) \rightarrow  0 $ in $L_{loc}^2(\mathbb{R}^N)$ ) and, up to a subsequence,
$v_n^1 (\cdot + y_n^2) \rightharpoonup w^2 $ in $H^1(\mathbb{R}^N )$  for some $w^2 \in B_m \backslash \{0\}.$
Noting that
$u_n (\cdot + y_n^2) =v_n ^1(\cdot + y_n^2)+w^1 (\cdot-y_n^1 + y_n^2)   \rightharpoonup w^2$ in $H^1(\mathbb{R}^N )$ ,
by (4.2) and arguing as above, we deduce that $P(w^2) = 0$ and thus $I(w^2) > 0.$
 Set $ v_n^2 =v_n^1-w^2(\cdot - y_n^2) =u_n- \Sigma_{l=1}^2 w^l( \cdot - y_n^l).$
Consequently,
$$  \lim_{n \rightarrow \infty}  \| \nabla v_n^2 \|_2^2= \lim_{n \rightarrow \infty}  \| \nabla u_n\|_2^2 - \Sigma_{l=1}^2  \|\nabla  w^l( \cdot - y_n^l) \|_2^2 ,$$
and
$$ 0>\lim_{n \rightarrow \infty}   I( v_n^1)=   I( w^2)+ \lim_{n \rightarrow \infty}   I( v_n^2) > \lim_{n \rightarrow \infty}   I( v_n^2).$$
We can obtain an infinite sequence $\{w^k\} \subset  B_m \backslash \{0\}$ such that
$P(w^k ) = 0$ and
$$\Sigma_{l=1}^k  \|\nabla w^l \|_2^2 \leq \lim_{n \rightarrow \infty}\| \nabla u_n\|_2^2 <+\infty,$$
for every $k \in \mathbb{N}^+$.
 However, this is in contradiction with $ \|\nabla w^l  \|_2 \geq \delta .$
 Therefore, the claim follows.

Let $s := \|w^1\|_2^2\in (0,m].$
Note that  $ \lim_{n \rightarrow \infty}  I(v_n^1) \geq0$ and  $w^1 \in  \mathcal{P}_s .$
 By (4.5), we have
$$E_m = I(w^1) + \lim_{n \rightarrow \infty}  I(v_n^1) \geq I(w^1)\geq E_s .$$
Since $E_m$ is nonincreasing in $m > 0$ due to Lemma 3.2, we have
$$I(w^1) = E_s = E_m, \eqno(4.6)$$
and
$$\lim_{n \rightarrow \infty}  I(v_n^1) = 0. \eqno(4.7)$$
  (4.3), (4.6) and Lemma 3.4 give $\mu \geq 0.$

To prove  $s = m,$ let us
prove that $\mu$ is positive.

The condition $(H3)$ gives that
$(N+\alpha)F(t) -(N-2)  f (t)t > 0$ for any $t \neq 0.$
 Thus, by
Pohozaev identity corresponding to (4.3), we obtain
$$\mu=\frac{ 1}{2s}\int_{\mathbb{R}^N} (I_\alpha \ast F(w^1)) [(N+\alpha)F(w^1)-(N-2)f(w^1)w^1]dx > 0. \eqno(4.8)$$

 Hence, we can obtain that $ \mu> 0.$

 If $s < m, $ from (4.3), (4.6) and
Lemma 3.4, it follows that
$I(w^1) = E_s > E_m$
which is in contradiction with (4.6).
Therefore, $s := \| w^1\|_2^2= m $ and then $\|v_n^1\|_2^2\rightarrow 0$ via (4.4).
Note that $\lim_{n \rightarrow \infty}\int_{\mathbb{R}^N} (I_\alpha \ast F(v_n^1))F(v_n^1)dx=0 $ due to Lemma 2.2 (ii ).
We observe that
$\|\nabla v_n^1\|_2^2\rightarrow 0$ due to  (4.7) and thus
$u_n(\cdot + y_n^1 )  \rightarrow w^1$ strongly in $H^1(\mathbb{R}^N).$
The proof of the lemma is completed.
  \hfill$\Box$

{\bf   Proof of Theorem 1.1.} From Lemmas 4.2 and 2.7 (iv), it follows that there exists a bounded Palais-Smale sequence
$\{u_n\} \subset  \mathcal{P}_m$ for the constrained functional $I|_{S_m}$ at the level $E_m > 0.$
Let $(H3)$ be  satisfied. Then Lemma 4.5 gives the existence of a
ground state $u \in S_m $ at the level $ E_m,$
and  $\|u_n^-\| _2^2 \rightarrow 0 $ thanks to Lemma 4.2.
By Lemma 4.5, there exists a nonnegative ground state
$u \in S_m$  at the level  $ E_m.$
Applying the strong maximum principle, $u > 0$ as required. \hfill$\Box$

{\bf Proof of Theorem 1.2.}
From Theorem 1.1, it follows that $E_m$ is reached by a ground state of $(1.1)$ with the associated
Lagrange multiplier being positive, and thus the function $m \mapsto  E_m$ is strictly decreasing on
$(0,\infty) $ thanks to  Lemma 3.4.  Combing with Lemmas 2.7, 3.1, 3.2, 3.5 and 3.6, the rest of the proof directly holds.
  \hfill$\Box$

\section{ Existence of  Radial solutions}

 Set $\sigma: H^1(\mathbb{R}^N  ) \rightarrow H^1(\mathbb{R}^N  )$ being the transformation
   $ \sigma(u) =-u $  and let $X\subset H^1(\mathbb{R}^N  ).$
   A set $A \subset X$ is called $\sigma$-invariant if $\sigma(A) = A.$
    A homotopy
$\eta : [0, 1] \times X \rightarrow X$ is $\sigma$- equivariant if
$\eta(t, \sigma(u)) = \sigma(\eta(t, u))$ for all $(t, u) \in [0, 1] \times X.$

{\bf Definition 5.1.}(\cite[Definition 7.1]{Ghoussoub}). Let $B$ be a closed  $\sigma$-invariant subset of $X\subset H^1(\mathbb{R}^N  ).$
 A class $\mathcal{G}$ of compact
subsets of $X$ is said to be a  $\sigma$-homotopy stable family with closed boundary B provided\\
(i ) every set in $\mathcal{G} $ is  $\sigma$-invariant,\\
(ii ) every set in $\mathcal{G}$ contains $B,$\\
(iii ) for any set $A \in \mathcal{G} $ and any  $\sigma$-equivariant homotopy $\eta \in C([0, 1]\times X, X)$ that satisfies
$\eta(t, u) = u$ for all $(t, u) \in (\{0\} \times X) \cup ([0, 1]\times B),$ one has $\eta(\{1\} \times A) \in \mathcal{G}.$

When $f$ is odd, we obtain that $s(-u) = s(u) $  due to Lemma 2.6 (iv)  and thus the constrained
functional
$$J (u) = I (s(u)\star u) = \frac{1}{2}e^{2s(u)}\int_{\mathbb{R}^N }|\nabla u|^2dx - \frac{1}{2}e^{-(N+\alpha)s(u) }
\int_{\mathbb{R}^N } \left(I_\alpha \ast F(e^{Ns(u)/2}u)\right)F(e^{Ns(u)/2}u)dx$$
is even in $u \in S_m. $
That is, $J$ is a $\sigma$-invariant functional on $S_m.$

{\bf Lemma 5.2.} {\it Let $\mathcal{G}$  be a $\sigma$-homotopy stable family of compact subsets of $S_m  \cap H_r^1$
 (with $B = \emptyset$) and set
$$E_{m,\mathcal{G}} := \inf_{A\in \mathcal{G}}\max_{u\in A}J (u).$$
If $E_{m,\mathcal{G}}> 0,$ then there exists a Palais-Smale sequence $\{u_n\} \subset \mathcal{P}_m  \cap H_r^1  $ for the
constrained functional $I|_{S_m  \cap H_r^1(\mathbb{R}^N )}$ at the level $E_{m,\mathcal{G}}.$}

 {\bf Proof.}
The proof is almost the same as Lemma 4.4.

Choose a sequence of $\sigma$-homotopy stable families of compact subsets of $S_m  \cap  H_r^1(\mathbb{R}^N )$
 (with $B = \emptyset$).
 Fix a sequence of finite dimensional linear subspaces $\{V_k\} \subset   H_r^1(\mathbb{R}^N )$
such that $V_k\subset V_{k+1},$ $ \dim V_k = k$ and $\bigcup  _{k\geq 1}V_k$  is dense in $H_r^1(\mathbb{R}^N ),$
 and set $\pi_k$ being the
orthogonal projection from $ H_r^1(\mathbb{R}^N )$
  onto $V_k .$

{\bf Definition 5.3. } {\it For any nonempty closed $\sigma$-invariant set $A \subset H^1(\mathbb{R}^N ),$  the genus of A is
defined by
$$\text{Ind}(A) := \min\left\{k \in \mathbb{N}^+| \ \  \exists  \  \varphi : A \rightarrow \mathbb{ R}^k  \backslash \{0\}, \varphi \text{ is odd and continuous}\right\}.$$
}
\qquad Let Ind$(A) =\infty$ if such $ \varphi$ does not exist,
and set Ind$(A) = 0 $ if  $ A = \emptyset.$
Denote $\Sigma$ by the family of compact $\sigma$-invariant subsets of $S_m  \cap  H_r^1(\mathbb{R}^N )$.
For each $k \in \mathbb{N}^+$, set
$$\mathcal{G}_k:= \{A \in \Sigma | \  \text{Ind}(A) \geq k\},$$
and
$$E_{m,k} := \inf_{A\in \mathcal{G}_k}\max_{u\in A}J (u).$$

{\bf Lemma 5.4.}  {\it  (i ) For any $k \in \mathbb{N}^+,$
$$\mathcal{G}_k \neq \emptyset,$$
and $ \mathcal{G}_k $ is a $\sigma$-homotopy stable family of compact subsets of $S_m  \cap  H_r^1(\mathbb{R}^N )$
 (with $B = \emptyset$). \\
(ii ) $E_{m,k+1} \geq E_{m,k} > 0$ for any  $k \in \mathbb{N}^+.$
}

{\bf Proof.}
  (i ) For each  $k \in \mathbb{N}^+,$
$S_m \cap V_k \subset \Sigma,$ combing with the basic properties of the genus, we have
$Ind(S_m \cap V_k) = k$
and thus $ \mathcal{G}_k \neq \emptyset$ .
  The rest follows from  Definition 5.1 and   again the basic properties of the
genus.\\
(ii )  $ E_{m,k}$ is well defined due to Item (i ).
For any $A\in  \mathcal{G}_k ,$ since $s(u)\star u \in \mathcal{P}_m$
for all $u \in A$ , it follows  from Lemma 2.6 (iii ) and Lemma 2.7 (iii ) that
$$  \max_{u\in A}J (u)= \max_{u\in A}I (s(u)\star u) \geq \inf_{v\in \mathcal{P}_m}
I (v) > 0,$$
and thus $ E_{m,k} > 0 .$
 Noting that $\mathcal{G}_{k+1} \subset \mathcal{G}_k,$  we have $E_{m,k+1}\geq  E_{m,k} .$
\hfill$\Box$

{\bf Lemma  5.5.} {\it Suppose that $\{u_n \} \subset H_r^1(\mathbb{R}^N).$
If $u_n \rightharpoonup u$ in $H_r^1(\mathbb{R}^N),$
then
$$     \lim _{n \rightarrow +\infty} \int_{\mathbb{R}^N}(I_\alpha \ast F(u_n))f(u_n)u_ndx=   \int_{\mathbb{R}^N} (I_\alpha \ast F(u)) f(u) udx.$$
}
\quad{\bf Proof.} Note that for any $ \varepsilon > 0,$ there exists  $C_\varepsilon>0$ such that
 $  |f(t)| \leq C_\varepsilon |t|^{ \frac{\alpha+2}{N}}+\varepsilon |t|^{{\frac{2+\alpha}{N-2}}}$
  for all $t \in \mathbb{R}.$
By Sobolev embedding theorem, $u_n \rightarrow u$ in $L^r(\mathbb{R}^N )$ ($2<r<2^*$), then
\begin{align*}
& \int_{\mathbb{R}^N} | f(u_n)  (u_n-u)|^{\frac{ 2N}{N+\alpha}}dx\\
 \leq &C C_\varepsilon [ \int_{\mathbb{R}^N}  |u_n| ^{ \frac{ \alpha +2}{N}  \frac{ 2N}{N+\alpha} \frac{ N+\alpha+2}{\alpha+2}}dx]^  \frac{\alpha+2}{ N+\alpha+2}
[ \int_{\mathbb{R}^N}    |u_n-u| ^{\frac{ 2N}{N+\alpha}\frac{ N+\alpha+2}{N}}dx]^\frac{N}{ N+\alpha+2}\\
+& C \varepsilon[ \int_{\mathbb{R}^N}  |u_n| ^{ \frac{ \alpha +2}{N-2}  \frac{ 2N}{N+\alpha} \frac{ N+\alpha}{\alpha+2}}dx]^  \frac{\alpha+2}{ N+\alpha}
[ \int_{\mathbb{R}^N}    |u_n-u| ^{\frac{ 2N}{N+\alpha}\frac{ N+\alpha}{N-2}}dx]^\frac{N-2}{ N+\alpha}\\
\rightarrow &0.
\end{align*}
Since $ |f(u_n)|^{\frac{2N}{N+\alpha}}  \in L^{{\frac{N+\alpha}{\alpha+2}}}( \mathbb{R}^N)$
and $ f(u_n) \rightarrow f(u) $ a.e. in $\mathbb{R}^N,$
we have $ | f(u_n)-f(u)|^{\frac{2N}{N+\alpha}}\rightharpoonup 0$ in $ L^{{\frac{N+\alpha}{\alpha+2}}}( \mathbb{R}^N).$  By the fact that $| u|^{\frac{2N}{N+\alpha}}$ in $L^{\frac{N+\alpha}{N-2}}(\mathbb{R}^N ),$
 we obtain
$$\int_{\mathbb{R}^N} |(f(u_n)-f(u)) u|^{\frac{2N}{N+\alpha}}dx  \rightarrow  0 .$$
Hence,
\begin{align*}
\int_{\mathbb{R}^N} | f(u_n)u_n -f(u)u|^{\frac{2N}{N+\alpha}}dx
\leq & C\int_{\mathbb{R}^N} |( f(u_n) -f(u) )u|^{\frac{ 2N}{N+\alpha}}dx
+C\int_{\mathbb{R}^N} | f(u_n)  (u_n-u) |^{\frac{ 2N}{N+\alpha}}dx\\
\rightarrow &0.
\end{align*}
Noting that $I_\alpha $ is a linear continuous map from $L^{\frac{2N}{N+\alpha}}(\mathbb{R}^N)$ to $L^{\frac{2N}{N-\alpha}}(\mathbb{R}^N),$ and
$F(u_n)$ is bounded in $L^{\frac{2N}{N+\alpha}}(\mathbb{R}^N)$,
we have $I_\alpha \ast F(u_n) \rightharpoonup I_\alpha \ast F(u ) $ in $L^{\frac{2N}{N-\alpha}}(\mathbb{R}^N).$

Consequently, we have
\begin{align*}
 \lim_{n \rightarrow \infty}\int_{\mathbb{R}^N} &[(I_\alpha \ast F(u_n))f(u_n)u_n-(I_\alpha \ast F(u))f(u)u] dx\\
\leq  &\lim_{n \rightarrow \infty}\int_{\mathbb{R}^N} |(I_\alpha \ast F(u_n))[f(u_n)u_n-f(u)u]| dx\\
&+ \lim_{n \rightarrow \infty}\int_{\mathbb{R}^N}| [(I_\alpha \ast F(u_n))-(I_\alpha \ast F(u))]f(u)u|dx\\
\leq & \lim_{n \rightarrow \infty} C\| F(u_n)\|_{\frac{2N}{N+\alpha}} \|f(u_n)u_n-f(u)u\|_{\frac{2N}{N+\alpha}}
+0\\
= & 0.
\end{align*}
We have the assertion.
\hfill$\Box$

{\bf Lemma  5.6.} {\it Assume that  $\{u_n\} \subset  S_m \cap  H_r^1(\mathbb{R}^N) $ be any bounded Palais-Smale sequence for the constrained functional $I|_{ S_m \cap  H_r^1(\mathbb{R}^N)}$  at an arbitrary level $c > 0 $ such that $ P(u_n) \rightarrow  0. $ Then there exists  $u  \in S_m \cap  H_r^1(\mathbb{R}^N) $
and $\mu > 0 $ such that, up to the extraction of a subsequence,
$u_n \rightarrow u $ strongly in $ H_r^1(\mathbb{R}^N) $
and $- \triangle u + \mu u =   (I_\alpha  \ast F(u))f (u).$}

{\bf Proof.}
 Recall that the sequence $\{u_n\}$ is bounded in $ H_r^1(\mathbb{R}^N) $.
  Up to a subsequence, there exists  $u \in H_r^1(\mathbb{R}^N) $ such that $u_n \rightharpoonup u $ in $ H_r^1(\mathbb{R}^N) $,
  $ u_n \rightarrow  u $ in $L^p(\mathbb{R}^N)$ for any $p \in (2, 2^*)$,
  and $u_n \rightarrow  u$ almost everywhere in $\mathbb{R}^N.$
  Noting that $\| d I(u_n)\|_{u_n ,*} \rightarrow 0,$   from \cite[Lemma 3]{Berestycki2}, it follows that
$$- \triangle u_n+ \mu_n u_n -   (I_\alpha  \ast F(u_n))f (u_n) \rightarrow 0 \ \text{in} \   (H_r^1(\mathbb{R}^N))^* ,\eqno(5.1)$$
where
$$\mu_n := \frac{1}{m}\left( \int_{\mathbb{R}^N}  (I_\alpha  \ast F(u_n))f (u_n)u_n dx-\|\nabla u_n\|_2^2\right).$$
Without loss of generality, one may assume that $\mu_n  \rightarrow \mu$  for some $\mu \in \mathbb{R}.$
 A similar argument to the proof of (4.3) and  the Palais principle of symmetric criticality \cite{Palais} show that
$- \triangle u + \mu u =   (I_\alpha  \ast F(u))f (u).$

Claim: $u \neq 0.$

Assume by contradiction that $u = 0 ,$
 then $u_n \rightarrow 0$  in $ L^{2(N+\alpha+2)/(N+\alpha)} (\mathbb{R}^N ).$
By Lemma 2.2 (ii ) and that $P(u_n) \rightarrow 0,$
 we have $  \int_{\mathbb{R}^N}  (I_\alpha  \ast F(u_n))F (u_n) dx\rightarrow 0$  and
$$  \|\nabla u_n\|_2^2= P(u_n)+ \frac{N}{2} \int_{\mathbb{R}^N}  (I_\alpha  \ast F(u_n))\widetilde{F} (u_n)dx  \rightarrow 0.$$
This gives
$$ c= \lim_{n \rightarrow \infty }I (u_n)=\frac{1}{2}\lim_{n \rightarrow \infty }\|\nabla u_n\|_2^2-\frac{1}{2}\lim_{n \rightarrow \infty } \int_{\mathbb{R}^N}  (I_\alpha  \ast F(u_n))F (u_n)dx =0 $$
which   is in contradiction with $c > 0.$
Since $u \neq 0 ,$  by a similar argument to the proof of  (4.8),  we have
$$\mu  := \frac{1}{2\|u\|_2^2}\left(  \int_{\mathbb{R}^N} (I_\alpha  \ast F(u ))[(N+\alpha) F(u )-(N-2)f (u )u ]dx\right)>0.$$
Noting that  $u_n \rightharpoonup u$ in $H_r^1(\mathbb{R}^N)$,
by Lemma 5.5 and (5.1),  we have
\begin{align*}
 \int_{\mathbb{R}^N}  | \nabla u |^2dx+\mu \int_{\mathbb{R}^N}  |u |^2dx
 =&  \int_{\mathbb{R}^N} (I_\alpha  \ast F(u ))f(u)udx\\
 =&   \lim_{n \rightarrow \infty} \int_{\mathbb{R}^N} (I_\alpha  \ast F(u_n ))f(u_n)u_ndx\\
 =&\lim_{n \rightarrow \infty} \int_{\mathbb{R}^N}  | \nabla u_n |^2dx+ \mu \lim_{n \rightarrow \infty}\int_{\mathbb{R}^N}  |u_n |^2dx.
\end{align*}
Because $\mu > 0, $ we obtain
$$\lim_{n \rightarrow \infty} \int_{\mathbb{R}^N}  | \nabla u_n |^2dx=  \int_{\mathbb{R}^N}  | \nabla u |^2dx,
\ \  \lim_{n \rightarrow \infty}\int_{\mathbb{R}^N}  |  u_n |^2dx =m = \int_{\mathbb{R}^N}  |u |^2dx,$$
and then $u_n \rightarrow   u$ in $H_r^1(\mathbb{R}^N).$
\hfill$\Box$

{\bf Lemma 5.7.} {\it $ E_{m,k} \rightarrow +\infty $ as $k \rightarrow +\infty.$}

{\bf Lemma 5.8.}  {\it  For every $c > 0,$ there exists $\rho = \rho(c) > 0$ small enough and $k(c) \in  \mathbb{N}^+$
sufficiently large such that for any $k \geq k(c) $  and any $u \in \mathcal{P}_m \cap H_r^1(\mathbb{R}^N ) $
 satisfying
$$I (u) \geq c \ \text{if} \   \| \pi_k u\|_{H^1 }\leq \rho.$$}
\quad{\bf Proof.} Assume by contradiction  that there exists $c_0 > 0  $ such that for any  $\rho  > 0$ and any
$k \rightarrow +\infty ,$ there exist  $l = l(\rho, k)\geq k$ and $u = u(\rho, k) \in \mathcal{P}_m \cap H_r^1(\mathbb{R}^N )$  such that
$ \| \pi_l u\|_{H^1 }\leq \rho,$ but $ I (u) < c_0.$
Thus, there exists a strictly increasing sequence $\{k_j\} \subset \mathbb{N}^+$
  (and
$\lim_{ j\rightarrow \infty}k_j =\infty$)  and a sequence $ \{u_j\} \subset \mathcal{P}_m \cap H_r^1(\mathbb{R}^N )$
such that
$ \| \pi_{k_j} u_j\|_{H^1 }\leq \rho$
and $ I (u_j) < c_0 $ for any $ j \in \mathbb{N}^+$.
Noting that $\{u_ j \}$ is bounded in $H_r^1(\mathbb{R}^N )$
  due to Lemma 2.7 (iv), up to a subsequence,
there exists $u\in H_r^1(\mathbb{R}^N )$
  such that
$u_j \rightharpoonup u$ in $H_r^1(\mathbb{R}^N )$
  and $u_j \rightharpoonup u$ in $L^2(\mathbb{R}^N  ).$

Claim that $u = 0.$

 Indeed, as $k_j\rightarrow \infty, $ we have
$\pi_{k_j} u \rightarrow  u$ in $L^2(\mathbb{R}^N  )$ and thus
$$\langle \pi_{k_j} u_j , u\rangle_ {L^2(\mathbb{R}^N  )}
=\langle   u_j ,\pi_{k_j} u\rangle_{L^2(\mathbb{R}^N  )}
\rightarrow \langle u,u \rangle _{L^2(\mathbb{R}^N  )}
\ \text{ as } j\rightarrow \infty.$$
Noting that $\pi_{k_j} u_j \rightarrow  0 $ in $ L^2(\mathbb{R}^N  )$, we  have
$\|u\| _ 2
= \lim_{ j\rightarrow \infty}\langle \pi_{k_j} u_j , u\rangle_ {L^2(\mathbb{R}^N  )}
= 0,$
which implies that  the claim follows.

Now, up to a subsequence, $ u_j  \rightarrow  0$  in $L^{{\frac{2(N+\alpha+2)}{N+\alpha}}} (\mathbb{R}^N )  $
thanks to  the compact inclusion.
Since $ \{u_j\} \subset \mathcal{P}_m \cap H_r^1(\mathbb{R}^N ),$
 by Lemma 2.2 (ii ),
we obtain
$$\int_{\mathbb{R}^N} |\nabla u_ j |^2dx = \frac{N}{2}\int_{\mathbb{R}^N}(I_\alpha \ast F(u_j)) \widetilde{F}(u_j)
dx \rightarrow 0\ \text{ as} \ j \rightarrow \infty,$$
which   is in contradiction with  Lemma 2.7 (ii ).
The conclusion holds.
\hfill$\Box$

{\bf Proof of Lemma 5.7.} We show by contradiction that
$$ \liminf_{k \rightarrow \infty}E_{m,k} < c \ \text{for some} \ c > 0. \eqno(5.2)$$
Set  $\rho(c) > 0$ and $k(c) \in \mathbb{N}^+$ being the numbers given by Lemma 5.8.
 By (5.2),
there exists $k > k(c)$ such that $E_{m,k} < c. $
We can then find $A \in  \mathcal{G}_k$
(that is $A \in \Sigma   $ and Ind$(A) \geq k$) such that
$\max_{u\in A}I (s(u)\star u) =\max_{u\in A}J (u) < c $ thanks to the definition of $E_{m,k}. $
Note that Lemma 2.6 (iii ) and (iv) give  the mapping $\varphi  : A \rightarrow  \mathcal{P}_m \cap H_r^1(\mathbb{R}^N )$ defined by
$\varphi(u) = s(u)\star u$ is odd and continuous.
Then we obtain that  $  \varphi(A) \subset \mathcal{P}_m \cap H_r^1(\mathbb{R}^N ),$
 $ \max_{v\in \varphi(A) }I (v) <c $  and
$$\text{Ind}( \overline{\varphi (A)} )\geq \text{Ind}(A) \geq k > k(c). \eqno(5.3)$$
By Lemma 5.8,  we have $ \inf_{v\in \overline{A}}\|\pi_{k(c)} v\|_{H^1} \geq \rho(c) > 0. $
Let
$\psi(v) = \frac{1}{\|\pi_{k(c)} v\|_{H^1} }\pi_{k(c)} v$
  for any $v \in A,$
we obtain an odd continuous mapping $\psi : A \rightarrow \psi(\overline{A})\subset V_{k(c)} \backslash \{0\} $ and thus
Ind$(\overline{A}) \leq  $Ind$(\psi(\overline{A})) \leq k(c)$
which is in contradiction with (5.3).
Then,   $E_{m,k} \rightarrow +\infty$ as $k \rightarrow +\infty.$
\hfill$\Box$

{\bf Proof of Theorem 1.3.}
\quad For each $k \in \mathbb{ N}^+,$ from Lemmas 5.2 and 5.4, it follows that there exists a Palais-Smale sequence
 $\{u_n^k\}_{n=1}^\infty \subset \mathcal{P}_m \cap H_r^1(\mathbb{R}^N )$
  of the constrained functional  $I|_{S_m \cap H_r^1(\mathbb{R}^N )} $ at the level $E_{m,k} > 0. $
The sequence is bounded in $H_r^1(\mathbb{R}^N ) $ due to Lemma 2.7 (iv).
Then, we see that by  Lemma 5.6  that $(1.1)$ has a radial solution $u_k$ with $I (u_k ) = E_{m,k} . $
 By Lemma 5.4 (ii ) and Lemma 5.7,  we then have
$I (u_{k+1}) \geq I (u_k) >0 $ for  any  $k \geq 1$
and $I (u_k ) \rightarrow+\infty.$
\hfill$\Box$

\smallskip
{\small
\noindent {\bf Acknowledgements.}
S.W. Ma was supported by National Natural Science Foundation of China
(Grant Nos.11571187, 11771182)

\end{document}